\documentclass[11pt]{article}
\usepackage{graphicx}
\usepackage{amsfonts}
\newtheorem{algorithm}{Algorithm}[section]
\newtheorem{remark}{{Remark}}[section]
\def\cir{\makebox[0.5cm]{$\circ$}}
\def\ind{\hspace{0.25in}}

\def\fP{\em \sf}
\def\IR{\mathop{\mbox{\rm Iround}}\nolimits}
\def\real{\mathbb{R}}

\newcommand{\beq}{\begin{eqnarray}}
\newcommand{\eeq}{\end{eqnarray}}
\newcommand{\bd}{\begin{description}}
\newcommand{\ed}{\end{description}}

\newcommand{\realtre}{\mathbb R^3}
\newcommand{\intrtre}{\int_{\realtre}}

\newcommand{\eps}{\varepsilon}

\newcommand{\be}{\begin{equation}}
\newcommand{\ee}{\end{equation}}
\newcommand{\ba}{\begin{array}}
\newcommand{\ea}{\end{array}}

\renewcommand{\Sigma}{{\bar{\sigma}}}
\title{Adaptive and Recursive Time
Relaxed Monte Carlo methods for rarefied gas dynamics}

\author{Stefano Trazzi
\thanks{University of Ferrara, Department of Mathematics, Via
Machiavelli 35, 44100 Ferrara Italy. E-mail: {\tt
trazzi@dm.unife.it }},\and Lorenzo Pareschi
\thanks{University of Ferrara, Department of Mathematics,
Via Machiavelli 35, 44100 Ferrara Italy. E-mail: {\tt
lorenzo.pareschi@unife.it}},\and Bernt Wennberg \thanks{Chalmers
University of Technology, Department of Mathematics, SE-412 96
G\"oteborg, Sweden. E-mail: {\tt wennberg@chalmers.se}}}

\begin{document}

\date{October 3, 2008}

\maketitle

\begin{abstract}
Recently a new class of Monte Carlo methods, called Time Relaxed
Monte Carlo (TRMC), designed for the simulation of the Boltzmann
equation close to fluid regimes have been introduced \cite{PR-r}.
A generalized Wild sum expansion of the solution is at the basis
of the simulation schemes. After a splitting of the equation the
time discretization of the collision step is obtained from the
Wild sum expansion of the solution by replacing high order terms
in the expansion with the equilibrium Maxwellian distribution; in
this way speed up of the methods close to fluid regimes is
obtained by efficiently thermalizing particles close to the
equilibrium state. In this work we present an improvement of such
methods which allows to obtain an effective uniform accuracy in
time without any restriction on the time step and subsequent
increase of the computational cost. The main ingredient of the new
algorithms is recursivity \cite{Pareschi-Wennberg-r}. Several
techniques can be used to truncate the recursive trees generated
by the schemes without deteriorating the accuracy of the numerical
solution. Techniques based on adaptive strategies are presented.
Numerical results emphasize the gain of efficiency of the present
simulation schemes with respect to standard DSMC methods.

\end{abstract}

%\tableofcontents


{\bf AMS} 65C05, 76P05, 82C80

% ----------------------------------------------------------------
\section{Introduction}
Computations for rarefied gas dynamics (RGD) in engineering
applications are most frequently performed using Monte Carlo
methods. The Bird's method has been particularly successful for a
wide range of applications \cite{Bird3-r}. One of the major
drawbacks of these methods is the difficulty to compute the
simulation of rarefied gases that are close to the fluid dynamic
limit, since in such regime the collisional time becomes very
small. A nondimensional measure of the significance of collisions
is given by the Knudsen number $\varepsilon$, which is small in
the fluid dynamic limit and large in the rarefied state. For small
Knudsen numbers most Monte Carlo methods lose their efficiency
because they are forced to operate on a very short time scale. The
aim of this paper is to introduce a new Monte Carlo method that is
robust in the fluid dynamic limit, by which we mean that it is
accurate and efficient for a full range of Knudsen numbers.
Alternative approaches to the derivation of efficient Monte Carlo
methods for the simulation of the Boltzmann equation have been
presented by several authors (see for example \cite{BH-r, BS,
HashHassan-r, HH-r, LiZa-r, CPmc-r, SB-r} and the references
therein). Most of these approaches focus on the variance reduction
of the method using information coming from the macroscopic scale
(equilibrium states) at different levels.

The starting point of the construction of this new family of
schemes is the time relaxed discretization of the Boltzmann
equation by the Wild sum \cite{Wild-r}. Given a set of particles
one tries to split it into subsets of particles according to the
probabilities given by the coefficients in the sum. The remaining
particles are sampled from a Maxwellian correspondent to the local
equilibrium. This can be done with a recursive algorithm in an
efficient way. As a result the method does not contain any time
discretization error (similarly to Bird's method) on the contrary
to Nanbu-Babovsky method or classical TRMC methods
\cite{Nanbu80-r, PR-r}.

The goal of recursive TRMC (TRMC-R) methods is to construct simple
and efficient numerical methods for the solution of the Boltzmann
equation in regions with a large variation in the mean free path.
As a consequence the TRMC-R methods have the following features:

\begin{itemize}
\item for large Knudsen numbers, the TRMC-R methods behave as a Bird's method;

\item for intermediate Knudsen numbers the methods adapt the length of the collision trees in order
to speed up the computation time without degradation of accuracy;

\item in the limit of the very small Knudsen number, the collision step replaces the
distribution function by a local Maxwellian with the same moments.
The methods will behave as a stochastic kinetic scheme for the
underlying Euler equations of gas dynamics \cite{pullin78-r};

\item
mass, momentum, and energy are preserved.
\end{itemize}

The paper is divided into 5 sections. After this introduction we
will recall some basic facts about the Boltzmann equation and its
fluid-dynamic limit. Next in section 3 we will discuss the problem
of the time discretization of the Boltzmann equation and the basic
notations of TRMC method. Section 4 is devoted to a detailed
description of TRMC-R for Maxwell and Hard Sphere models. We
present also adaptive techniques to improve the efficiency without
degradation in accuracy. Finally in the last section we present
detailed numerical tests both for a space homogeneous case and for
a stationary shock wave, using the Hard Sphere Model. The results
show a marked improvement in the efficiency of computations given
by TRMC-R over Bird's method.

\section{The Boltzmann Equation}
The Boltzmann equation describes the evolution of a continuum of particles by mean of
three variables: the time $t$, the position $x$ and velocity $v$ of particles.

In this model, the density $f=f(x,v,t)$ of particles follows the equation

\begin{equation}
\frac{\partial f}{\partial t} + v \cdot \nabla_x f =
\frac{1}{\epsilon} Q(f,f),\quad  x\in \Omega \subset \real^3, v
\in \real^3,
\label{eq:BE-r}
\end{equation}

supplemented with the initial condition \be f(x,v,t=0) =
f_{0}(x,v). \label{eq:IC-r} \ee The function $f$ can depend on other
independent variables like an internal energy \cite{cercignani-r}.

In (\ref{eq:BE-r}) the parameter $\epsilon > 0$ is called
\emph{Knudsen number} and it is proportional to the mean free path
between collisions. The bilinear collisional operator $Q(f,f)$
which describes the binary collisions between particles in a
mono-atomic gas is given by

\begin{equation}
Q(f,f)(v)=\int_{R^3} \int_{S^2} B(f(v')f(v'_*)-f(v)f(v_*)) dv_* d \sigma,
\label{eq:Q-r}
\end{equation}
where for simplicity the dependence of $f$ on $x$ and $t$ has been
omitted.

In the previous expression $\sigma$ is a vector of the unitary
sphere $S^2 \subset R^3$. The collisional velocities $(v',v'_*)$
are associated to the velocities $(v,v_*)$ and to the parameter
$\sigma$ by the relations
\begin{equation}
v'=\frac{1}{2}(v+v_*+|q|\sigma),\quad
v'_*=\frac{1}{2}(v+v_*+|q|\sigma),
\end{equation}
where $q=v-v_*$ is the relative velocity.

The kernel $B$ is a non negative function which characterizes the
details of the binary interaction between particles. The classical
Variable Hard Spheres model used for hypersonic
flows in the upper-atmosphere is \[B(|q|,|q \cdot \sigma|)=K|q|^\alpha, \\
0 \leq \alpha < 1,\] where $K$ is a positive constant. The case
$\alpha=0$ corresponds to a \emph{Maxwellian gas}, while
$\alpha=1$ is called a \emph{Hard Sphere Gas}.

Since $f$ is a mass density in the phase space to obtain the
density $\rho=\rho(x,t)$ we have to integrate $f$ in $v$ \be
\rho=\intrtre fdv. \ee Similarly the gas velocity $u$ is
determined by \be \rho u =\intrtre fvdv, \ee and the gas
temperature by \be T=\frac{1}{3\rho}\intrtre(v-u)^2fdv. \ee
Finally we define \be E=\frac{1}{2}\intrtre|v|^2fdv = \frac32 \rho
T + \frac12 \rho u^2\ee the energy per volume density.

The collisional operator is such that the \emph{H-Theorem} holds
\be \intrtre Q(f,f)\log(f) dv \leq 0. \ee This condition implies
that each function $f$ in equilibrium (i.e. $Q(f,f)=0$) has
locally the form of a Maxwellian distribution \be M(\rho,u,T)(v) =
\frac{\rho}{(2\pi T)^{3/2}}\exp\left( - \frac{|u-v|^2}{2T}\right),
\label{eq:Max}\ee where $\rho,u,T$ are the density, the mean
velocity and the gas temperature.

Now, if we consider the Boltzmann equation (\ref{eq:BE-r}) and
multiply it for the elementary collisional invariants
${1,v,|v|^2}$ and integrate in $v$ we obtain \[\intrtre
Q(f,f)\phi(v) dv = 0,\qquad \phi(v)=1,v,|v|^2,\] which correspond
to conservation of mass, momentum and energy. If it is possible to
invert the order of derivation - integration, we get \be
\frac{\partial}{\partial t}\int \phi(v)dv +
\sum_{i=1}^{3}\frac{\partial}{
\partial x_i}\int v_i \phi fdv =0.
\ee

Unfortunately replacing $\phi$ by the functions $1,v,|v|^2$ we
obtain a differential equation system which is not close since it
involves higher order moments of the function $f(x,v,t)$. Formally
as $\eps \to 0$ the function $f$ is locally replaced by a
Maxwellian. In this case it is possible to compute $f$ from its
moments using (\ref{eq:Max}) thus obtaining to leading order the
closed Euler system of compressible gas dynamics

\be \frac{\partial \rho}{\partial t} +
\sum_{i=1}^{3}\frac{\partial}{\partial x_i} (\rho u_i)=0, \ee \be
\frac{\partial}{\partial t}(\rho u_j) +
\sum_{i=1}^{3}\frac{\partial} {\partial x_i}(\rho u_i u_j) +
\frac{\partial} {\partial x_j} p=0, \quad j=1,2,3 \ee \be
\frac{\partial E}{\partial t} + \sum_{i=1}^{3}
\frac{\partial}{\partial x_i}(Eu_i+pu_i)=0, \ee where $p=\rho T$.

\section{Time Relaxed schemes}
\subsection{Time discretizations}

The starting point is the usual first order splitting in time of
(\ref{eq:BE-r}), which consists of solving separately a purely
convective step (i.e., $Q\equiv 0$ in (\ref{eq:BE-r})) and a
collision step characterized by a space homogeneous Boltzmann
equation (i.e., $\nabla_x f \equiv 0$ in (\ref{eq:BE-r})). Clearly,
after this splitting, almost all the main difficulties are
contained in the collision step. For this reason, in what follows
we will fix our attention on the time discretization of the
homogeneous Boltzmann equation \be {\frac{\partial f}{\partial t}}
= \frac{1}{\varepsilon} Q(f, f) . \label{eq:HOM-r} \ee

As proposed in \cite{toscani-r}, a general idea for deriving
robust numerical schemes, that is, schemes that are
unconditionally stable and preserve the asymptotic of the
fluid-dynamic limit, for a nonlinear equation like
(\ref{eq:HOM-r}), is to replace high order terms of a suitable
well-posed power series expansion by the local equilibrium.  Here
we will briefly recall the schemes presented in \cite{toscani-r}.

\subsection{Derivation}
Let us consider a differential system of the type
\be
\frac{\partial f}{\partial t} = \frac{1}{\varepsilon} \left [ P(f,f) - \mu f \right ],
\label{eq:HOM1-r}
\ee
with the same initial condition (\ref{eq:IC-r}),
and where $\mu\not=0$ is a constant and $P$ a bilinear operator.

Let us replace the time variable $t$ and the function $f=f(v,t)$
using the equations
\be
\tau = (1 -e^ {-\mu t / \varepsilon} ), \qquad
F(v,\tau) = f(v,t) e^{\mu t / \varepsilon}.
\ee

Then $F$ is easily shown to satisfy

\be
\frac{\partial F}{\partial \tau} = \frac{1}{\mu} P(F,F)
\label{eq:HOM2-r}
\ee

with $F(v,\tau =0) = f_0(v)$.

Now, the solution to the Cauchy problem for (\ref{eq:HOM2-r}) can be sought in the form of a power series
\be
F(v,\tau) = \sum_{k=0}^\infty {\tau}^k f_{k}(v),\qquad
f_{k=0}(v) = f_0(v),
\ee

where the functions $f_{k}$ are given by the recurrence formula
\be
f_{k+1}(v) = \frac{1}{k+1}\sum_{h=0}^k \frac{1}{\mu} P ( f_{h},
f_{k-h}),\quad k=0,1,\ldots.
\label{eq:CF-r}
\ee

Making use of the original variables, we obtain the following formal
re\-presentation of the solution to the Cauchy problem (\ref{eq:HOM-r}):
\be
f(v,t) = e^{{-\mu t / \varepsilon } }\sum_{k=0}^\infty
( {1-e^{-\mu t / \varepsilon}} )^k f_k(v) .
\label{eq:WS-r}
\ee

The method was originally developed by Wild \cite{Wild-r,CCG-r} to solve the
 Boltzmann equation for Maxwellian molecules.

From this representation, a class of numerical schemes can be naturally
derived.

In \cite{toscani-r}, the following class of numerical schemes, based on a
suitable truncation for $m \geq 1$ of (\ref{eq:WS-r}), has been constructed:
\be
f^{n+1}(v) = e^{{-\mu \Delta t / \varepsilon } }\sum_{k=0}^m
( {1-e^{-\mu \Delta t / \varepsilon}} )^k f_k^n(v) +
( {1-e^{-\mu \Delta t / \varepsilon}} )^{m+1} M(v),
\label{eq:WST-r}
\ee
where $f^n=f(n\Delta t)$ and $\Delta t$ is a small time interval.
The quantity $M$ (referred to as the local Maxwellian associated with $f$) is
the asymptotic stationary solution of the equation.

It can be shown that the schemes obtained in this way are of order
$m$ in time. Furthermore, these schemes satisfy the following
properties \cite{toscani-r}.

\begin{itemize}
\item[\rm  (i)] {\em Conservation.}

      If $P(f,g)$ is a nonnegative bilinear operator such
      that there exist some functions $\phi(v)$ with the following property,
      \be
         \int_{R^3} P(f,f)\phi(v)\, dv = \mu\int_{R^3}f\phi(v)\, dv,
      \ee
      and the initial condition $f^0$ is a nonnegative function, then
      $f^{n+1}$ is nonnegative for any $\mu\Delta t/\varepsilon$
      and satisfies
      \be
         \int_{R^3} f^{n+1}\phi(v)\, dv = \int_{R^3}f^n\phi(v)\, dv.
                                                        \label{eq:conserv-r}
      \ee
\item[\rm (ii)] {\em Asymptotic preservation (AP).}

      For any $m\ge 1$, we have
      \be
         \lim_{\mu\Delta t/\varepsilon \to\infty} f^{n+1} = M(v).
                                                        \label{eq:as_pres-r}
      \ee
\end{itemize}

In the case of the Boltzmann equation, with a collision kernel
bounded by $\bar\sigma$ taking $P(f,f)=Q(f,f)+\mu f$, with $\mu
\geq 4\pi\rho\bar\sigma$ the schemes guarantee the conservation of
mass, momentum and energy (by the first property) and the correct
solution near the fluid limit (i.e. $\epsilon \to 0$).

\section{TRMC-R methods}

%\label{sect5}
\subsection{Maxwellian case}

In order to simplify the derivation of the Recursive Time Relaxed
Schemes, we recall some basic facts on the algorithm for the
simple case of constant cross sections (Maxwellian molecules) as
proposed in \cite{PR-r}

First we note that in the case of Maxwell molecules the Wild sum
has a clear probabilistic interpretation. If $f(\Delta t)$ is the
velocity distribution of particle at time $\Delta t$, then taking
a particle at random from this distribution it might happen that
this particle has not collided one single time. The distribution
given this is just $f_0$ and the probability of this event is
$\exp({{-\mu \Delta t}/{\epsilon}})$. In the same way $f_1$ is the
velocity distribution for particles which have been involved in
exactly one collision, and the probability of that is
$(1-\exp({{-\mu \Delta t}/{\epsilon}}))(\exp({{-\mu \Delta
t}/{\epsilon}}))$. At least $f_m$ is the velocity distribution
given that exactly $m+1$ particles have been involved in their
collision history back to the initial time. To be able to find a
sample of $f_m$, we must assume that the densities $f_k, 0 \leq
m-1$ are all already known. Off course the only one of these that
is really known is $f_0$, the initial distribution.

A sample of $f_m$ can be determined in a recursive way. In order
to understand how the algorithm works, let's consider for example
the sampling from a starting density functions $f_k, \ k=1,2,3$.
From (\ref{eq:CF-r}) $f_1$ is obtained from the collisional
operator $P(f_0,f_0)$, $f_2$ from a combination of $P(f_0,f_1)$
and $P(f_1,f_0)$ with the same probability weight, where the terms
$f_1$ are constructed as seen before. Similarly $f_2$ is a
combination of $P(f_0,P(f_0,f_0))$ and $P(P(f_0,f_0),f_0)$.

In the same way it is possible to create the higher terms of the
Wild's Sum.
%For example the part of $f_3$ coming from $P(f_1,f_2)$
%will be created by $P(P(f_0,f_0),P(f_0,P(f_0,f_0)))$ and
%$P(P(f_0,f_0),P(P(f_0,f_0),f_0))$

A simple recursive Monte Carlo algorithm is the following\\
\begin{algorithm}[Recursive sampling for Maxwell molecules] ~\par
{\samepage\small
\begin{tabbing}
\ind 1. \= choose $n$ from a geometric distribution with parameter\\
                \> $ \tau = 1 - \exp(-\mu \Delta t/\epsilon)$\\
\ind 2. \> take a sample from the distribution with density $f_n$ \\
         \>        \= \cir \=  {\fP if} $n=0$, take $v$ from the initial density $f_0$  \\
         \>        \> \cir \>  {\fP else} proceed as follows \\
         \>        \>      \> -  \= choose $k \in \left\{0,1, \ldots, n-1\right\}$ with equal probability \\
         \>        \>      \>  - \>  take a sample $v_i$ from the density $f_k$\\
         \>        \>      \>      \> and $v_j$ from the density $f_{n-k-1}$ as in step 2\\
         \>       \>       \>  -   \>  perform the collision between $v_i$ and $v_j$, obtaining $v'_i$ and $v'_j$\\
         \>       \>       \> -    \>  then $v'_i$ and $v'_j$ are distributed according to the density $f_n$ \\

\end{tabbing}
}
\end{algorithm}

The post collisional velocities are computed through relations \be
    v'_i =  \frac{v_i+v_j}{2} + \frac{|v_i-v_j|}{2}\sigma,\quad
    v'_j =  \frac{v_i+v_j}{2} - \frac{|v_i-v_j|}{2}\sigma,
\label{eq:Maxwell_coll-r}
\ee
where $\omega$ is chosen uniformly in the unit sphere, according to
%\bd
%\item{2D:}
%\be
%   \omega = \left(\ba{c}\cos\theta \\ \sin\theta\ea \right), \quad
%       \theta = 2\pi\rand,

%                                                          \label{eq:n2d}
%\ee
%\item{3D:}
\be
   \sigma= \left(\ba{c}\cos\phi\sin\theta \\
                                     \sin\phi\sin\theta \\
                                     \cos\theta\ea\right), \quad
             \theta = \arccos(2\xi_1-1), \quad \phi = 2\pi\xi_2,
                                                          \label{eq:n3d-r}
\ee
%\ed
and $\xi_1,\xi_2$ are uniformly distributed random variables in ${[0,1]}$.
%\end{remark}

\begin{figure}[t]
\centering
\includegraphics[scale=1.0]{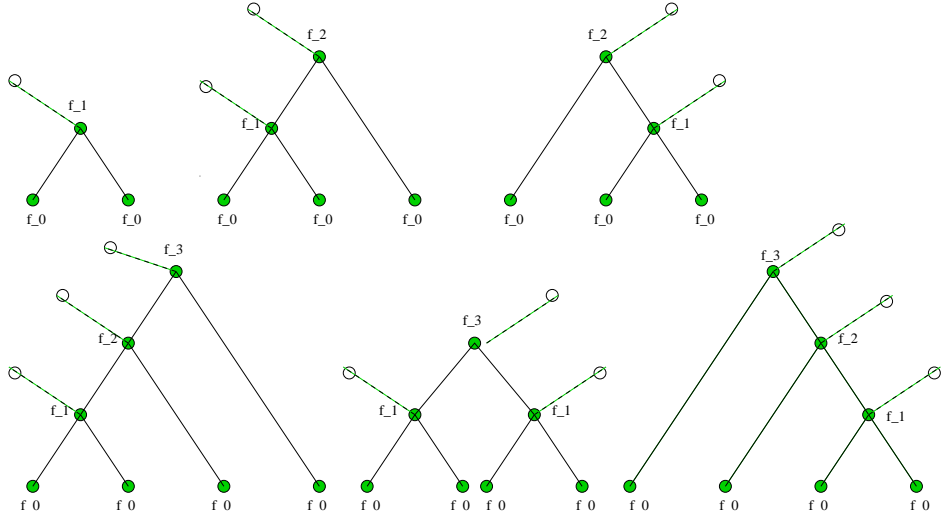}
\caption{Mc Kean graphs for $f_1$, $f_2$ and $f_3$}
\label{fig:a}
\end{figure}

It is useful to use a representation of the collision process
through the collision trees, sometimes called Mc Kean graphs
(Figure \ref{fig:a}).

It is now clear that the velocity distribution $f_k$ can be drawn
from the starting distribution $f_0$, choosing the different
collisional trees with the same probability in mean. Off course
two particles are produced in every collisional event, and only
one of these is used to complete the collisional process. So it is
natural to store the particle not used in a list, for its direct
utilization in a future collisional process. In this way we can
increase the efficiency of the method (if a particle sampled from
a velocity density $f_k$ already exist it will not be necessary to
obtain it by the complete collisional process) and guarantee the
exact conservation of moments.

In order to do this we first split a set of $N$ particles into $m$
collision sets $N_i$, $i=0,\ldots,m$ where each set $N_i$
characterizes the number of particles that will undergo $i$
collisions in a time step $\Delta t$. All particles having more
then $m$ collisions in the time step are thermalized (i.e.
replaced with an equilibrium particle taken from the local
Maxwellian), this number is denoted by $N_{m+1}$.
Accordingly to this we have the following algorithm where $N$ and $m$ are given.\\
\begin{algorithm}[Splitting particles into collision sets] ~\par
{\samepage\small
\begin{tabbing}
\ind 1. \= set $ \tau = 1 - \exp(-\mu \Delta t/\epsilon)$\\
        \> $n=0$, $\lambda_0 =1$, $\omega_0=(1-\tau)$, $\tilde N_0 =N $, $N_c = 0$\\
        \> $N_0 = \IR(\frac{\omega_0 \tilde N_0}{\lambda_0})$\\
\ind 2. \> repeat \\
         \>        \= - \=  $n=n+1$, \\
         \>        \> - \>  $\tilde N_n = \tilde N_{n-1} - N_{n-1}$ \\
         \>        \> - \> $\omega_n = (1-\tau)\tau^n $\\
         \>        \> - \> $\lambda_n = \lambda_{n-1} - \omega_{n-1}$\\
         \>        \> - \> $N_n = \IR(\frac{\omega_n \tilde N_n}{\lambda_n})$\\
         \>        \> - \> $N_c = N_c + N_n$\\
         \>until $\tilde N_n > 0$ and $n<m$ \\
\ind 3. \> if ($n=m$)\\
        \>         \> - \> $N_{m+1} = N-N_c-N_0$\\
        \> else\\
        \>          \> - \> $m=n$, $N_{m+1}=0$\\
        \> end if\\
\end{tabbing}

\label{alg:split}
}
\end{algorithm}

Here, by $\IR(x)$, we denote a suitable integer
rounding of a positive real number $x$. In our algorithm, we
choose
\[
   \IR(x) = \left\{\begin{array}{lll}
                     {[x]}     & \mbox{with probability} & {[x]}+1-x, \\
                     {[x]} + 1 & \mbox{with probability} & x-{[x]},
                   \end{array}
            \right.
\]
where $[x]$ denotes the integer part of $x$.

Note that, since we have a finite number of particles, a
consequence of the above splitting into collision sets is that the
maximum possible length of a collision process is fixed by the
initial number of particles $N$. The recursive collision algorithm
for Maxwell molecules where thermalization occurs accordingly to a
TR discretization of order $m$ is given here. To achieve exact
conservation of momentum and energy and a better computational
efficiency, the algorithm uses counters $c_i$ to keep track of the
number of particles stored in memory with a collision history of
$i$ collisions in a time step $\Delta t$ and not yet used in the
simulation.
\\
\begin{algorithm}[TRMC-R for Maxwell molecules] ~\par
{\samepage\small
\begin{tabbing}
\ind 1. \= compute the initial velocity of the particles, $ \left\{ v_i^0, i=1,\ldots,N  \right\}$\\
                \> by sampling them from the initial velocity $f_0$\\
\ind 2. \> split particles into collision sets as in algorithm (\ref{alg:split})\\
\ind 3. \> set counters $c_n = 0$ for $n=1,\ldots,m+1$\\
\ind 4. \> {\fP for} $n = m,\ldots,1$\\
        \>        \= take $N_n$ samples from the distribution with density $f_n$, according to\\
        \>        \> \cir  \= {\fP repeat}\\
        \>        \>       \>          \ind - \= choose $k \in \left\{0,1,\ldots, n-1\right\}$ with equal probability \\
        \>        \>       \>          \ind - \> if $k=0$ take $v_i$ from the initial density $f_0$ \\
        \>        \>       \>          \ind - \> else choose $v_i$ from the density $f_k$ \\
        \>        \>       \>            \>     \ind  \=if $c_k>0$ use a stored particle with a random choice\\
                \>        \>       \>            \>           \> set $c_k = c_k - 1$ and $N_k=N_k+1$\\
                \>        \>       \>            \>     \ind  \> else sample $v_i$ and $v_i^*$ from $f_k$ (recursively)\\
                \>        \>       \>            \>           \> $v_i^*$ is stored and then set $c_k=c_k+1$, $N_k=N_k-1$\\
        \>        \>       \>          \ind - \> if $n-k-1=0$ take $v_j$ from the initial density $f_0$ \\
        \>        \>       \>          \ind - \> else choose $v_i$ from the density $f_{n-k-1}$ \\
        \>        \>       \>            \>     \ind  \=if $c_{n-k-1}>0$ use a stored particle with a random choice\\
                \>        \>       \>          \>           \> set $c_{n-k-1} = c_{n-k-1} - 1$ and $N_{n-k-1}=N_{n-k-1}+1$\\
                \>        \>       \>          \>     \ind  \> else sample $v_i$ and $v_i^*$ from $f_{n-k-1}$ (recursively)\\
                \>        \>       \>          \>           \> $v_i^*$ is stored and then set $c_{n-k-1}=c_{n-k-1}+1$,\\                                \>        \>       \>          \>           \> $N_{n-k-1}=N_{n-k-1}-1$\\
        \>        \>       \>          \ind - \> perform the collision between $v_i$ and $v_j$ as in DSMC \\
        \>        \>       \>          \ind - \> $v'_i$ and $v'_j$ are random variables distributed \\
        \>        \>       \>                       \> according to the density $f_n$ \\
        \>        \>       \>          \ind - \> set $N_n=N_n - 2$ \\
        \>        \>                 \> {\fP until} $(N_n>0)$ \\
\ind    \> {\fP end for} \\
\ind 5. \> sample $N_{m+1}$ particles from the local Maxwellian\\

\end{tabbing}
\label{alg:m-r}
}
\end{algorithm}

\subsection{VHS collision kernels}
\label{sec:VHS}
The algorithm described for Maxwellian molecules can be extended to more general
collision kernels by using dummy collisions and acceptance-rejection technique.
This approach is equivalent to sample the post collisional velocity according to
$P(f,f)/ \mu $, where $\mu = 4\pi\Sigma$ and $\Sigma$ is an upper bound of the
scattering cross section for the given set of particles.

The upper bound $ \Sigma  $ should be chosen as small as
possible, to avoid inefficient rejection, and it should be
computed fast.

An optimal bound can be derived taking $\Sigma$ as \be \Sigma=
\max_{v_i,v_j} \sigma(|v_i-v_j|). \label{eq:SMAX-r} \ee However
this computation would be too expensive since it would require an
$O(N^2)$ operations. An upper bound of $\Sigma$ can be obtained by
taking  $\Sigma = \sigma(2\Delta v)$,
\[ \Delta v = \max_{i} |v_i-\bar{v}|,
\quad \bar{v}=\sum_i v_i/N. \]

In the VHS case the algorithm should be modified as follows\\
\begin{algorithm}[TRMC-R for VHS molecules] ~\par
{\samepage\small
\begin{tabbing}
\ind 1. \= compute the initial velocity of the particles, $ \left\{ v_i^0, i=1,\ldots,N  \right\}$\\
                \> by sampling them from the initial velocity $f_0$\\
\ind 2. \> split particles into collision sets as in algorithm (\ref{alg:split})\\
\ind 3. \> set counters $c_n = 0$ for $n=1,\ldots,m+1$\\
\ind 4. \> compute an upper bound $\Sigma$ of $\sigma_{ij} = \sigma(v_i,v_j)$\\
\ind 5. \> {\fP for} $n = m,\ldots,1$\\
        \>  \ind \= take $N_n$ samples from the distribution with density $f_n$, according to\\
        \>              \> point 4 of algorithm \ref{alg:m-r} where the dummy collision is performed as\\
        \>              \> in DSMC if $\Sigma \xi_1 < \sigma_{ij}$, with $\xi_1$ uniformly distributed random variable in ${[0,1]}$\\
\ind 6. \> sample $N_{m+1}$ particles from the local Maxwellian\\

\end{tabbing}
\label{alg:1-r}
}
\end{algorithm}

The algorithm is exactly conservative if combined with a suitable
scheme for sampling a set of particles with prescribed momentum
and energy from the Maxwellian, as proposed in \cite{pullin78-r}
or by the authors in \cite{partra-r}.

\begin{remark}{˜}
\begin{itemize}
\item
In the case of $m=1, 2$ the TRMC-R method corresponds to the first
and second order TRMC methods presented in \cite{PR-r} for the
simplest possible choice of the weights. Off course here we are
aiming at using much larger values of $m \gg 1$ for which a direct
extension of the algorithms presented in \cite{PR-r} is not
feasible. Larger values of $m$ produce higher accurate results and
then allow a larger time step when compared to \cite{PR-r,
partra-r}.

\item
Since the collision process changes the distribution function, it
is important to choose a correct upper bound $\Sigma$ of
$\sigma_{ij}$ in order to avoid discarding all the collisional
trees computed in case $\sigma_{ij}> \Sigma$. As an alternative
one can update the upper bound itself, after each collision as in
\cite{PR-r, partra-r}.
\end{itemize}
\end{remark}

\subsection{Adaptive technique for TRMC-R scheme }
In practical simulations the number $m$ can be very large,
depending on the Knudsen number and on the number of test
particles. Clearly small values of $m$ make the algorithm faster,
because the collision process is replaced by the projection to the
local Maxwellian equilibrium, but far from the fluid regime
keeping $m$ too small can produce less accurate results. In
practice, a maximum allowed value $m_{\max}$ of $m$ is fixed at
the beginning of the calculations; $m_{\max}$ represents the
maximum depth of a collision tree. The main problem is to choose
the right $m_{\max}$, in order to have the best combination
between efficiency and accuracy. The idea we develop is to use an
adaptive technique to choose the right maximum depths of the
collision trees, based on evaluating the distance of the solution
from the equilibrium through a suitable indicator. This can be
performed measuring the variation of some macroscopic variables
such as the fourth order moment or the components of the shear
stress tensor.

Let $S$ be the macroscopic variable selected accord\-ing to the particular physical problem.
Then define the quantity
\[
E_1 = \frac{|{S^{n+1,m_{\max}}-S^n}|}{|{S^n}|}
\]
that represents the relative variation at time step $n+1$ of the
macroscopic variable $S$ computed with the solution obtained using
$m_{\max}$ as maximum depth of the collision trees. If we fix an
interval $\left[\delta_1, \delta_2 \right]$, $0 < \delta_1 <
\delta_2$ then we can apply the following criteria in order to
accept or discard the solution at time step $n+1$

\begin{itemize}
\item if $E_1 < \delta_1$ the solution is accepted and $m_{\max}=m_{\max}/2$ in the next time step $n+2$;
\item if $\delta_1 \leq E_1 \leq \delta_2$ the solution is accepted and $m_{\max}$ is unchanged in the next time step $n+2$;
\item if $E_1 > \delta_2$ the solution is discarded and the calculation is performed again using $m_{\max}=2m_{\max}$.

\end{itemize}

Off course it is possible to use other techniques in order to
evaluate the distance from equilibrium \cite{Tiv-Rja}.

The algorithm works with optimal efficiency if the collisions
computed with the ``wrong" $m_{\max}$ are kept and reused with
$2m_{\max}$.

This can be done by starting with $m = m_{\min}$ and observing
that if
 \[\tilde{f}^{n+1,m} = \tau \sum_{k=0}^{m} (1-\tau)^k f_k + (1-\tau)^{m+1}M \] then
 \[{f}^{n+1,2m} = \tilde{f}^{n+1,m} + \tau \sum_{k=m+1}^{2m} (1-\tau)^k f_k + [(1-\tau)^{2m+1} - (1-\tau)^{m+1}]M \]

So if the test for the adaptive strategy impose to discard the solution we proceed as follows
\begin{itemize}
\item the collisions computed with  $m$ are kept;
\item the fraction $(1-\tau)^{m+1}M$ is discarded;
\item the fraction $\sum_{k=m+1}^{2m} (1-\tau)^k f_k$ is computed by the recursive collision
process;
\item the fraction $[(1-\tau)^{2m+1} - (1-\tau)^{m+1}]$ is sampled by a Maxwellian.
\end{itemize}

We observe that for estimation purposes the sampling from
Maxwellian can be substitute by the analytical computation of
local Maxwellian. This contributes to a better efficiency of the
adaptive method.

%\begin{center}
\begin{table}
\begin{center}
\begin{tabular}{|c|c|c|} \hline
           Definition  &  {left tree}    & {right tree} \\ \hline\hline
(\ref{def:1}) &  7                   &  7                \\
\hline (\ref{def:2}) &  3                   &  2
\\ \hline
(\ref{def:3}) &  3                   &  43/16            \\
\hline
\end{tabular}
\vskip .5cm \caption{Length of the collision trees}
\label{tab:len}
\end{center}
\end{table}
%\end{center}

\subsection{Truncation of the collision trees for TRMC-R scheme }

\begin{figure}[t]
\centering
\includegraphics[scale=0.45]{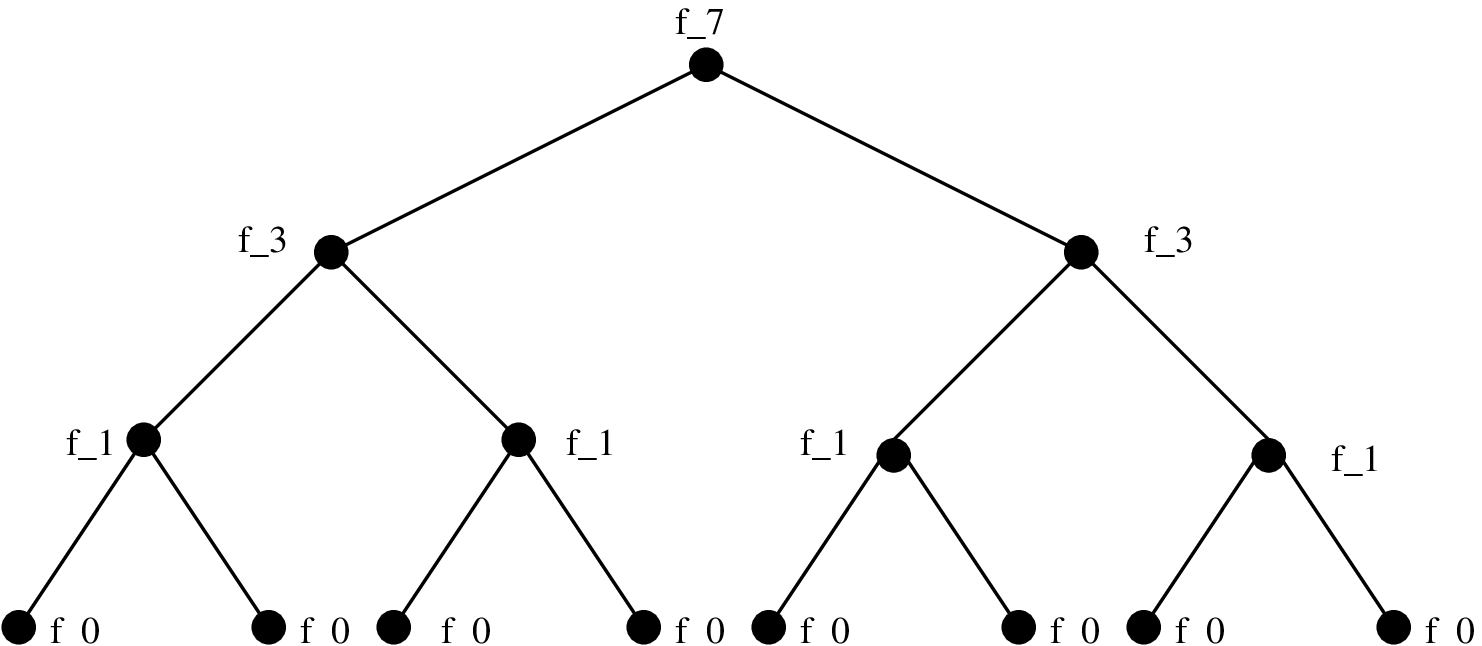}
\includegraphics[scale=0.65]{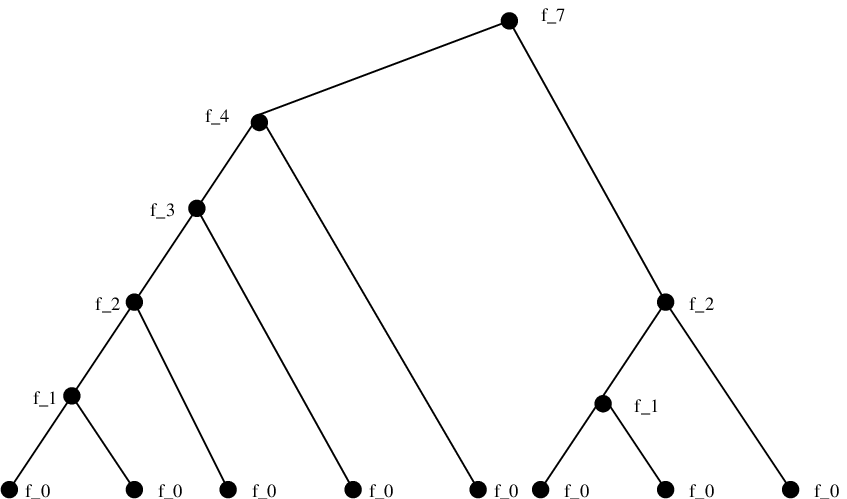}
\caption{Well balanced (left)  and not well balanced tree (right)}
\label{fig:b}
\end{figure}

Off course several definitions of the length of trees are
possible; the simplest one, that does not care about the shape of
the collision trees that generate a particle from the density
function $f_k$, is the one provided by \be L(k = h+j+1) = k.
\label{def:1}\ee The length corresponds to the coefficient of the
density function by which we sample the particle. The idea is to
sample directly from the local Maxwellian if, for the collision
process, $L(k)>m_{\max}$, $m_{\max}$ fixed, otherwise the whole
tree is kept and the collisions are performed. This simple
definition has been used into the algorithms described before.
Different definition of length $L$ can be done using recursivity
as \beq
    L(k = h+j+1) &=& 1 + \min\{L(h),L(j)\},\\
    \label{def:2}
    L(k = h+j+1) &=& 1 + \hbox{mean}\{L(h),L(j)\}.
    \label{def:3}
\eeq These two last definitions can be related to the concept of
``well balanced and not well balanced trees" (see Figure
\ref{fig:b}), i.e. if $L(k)>m_{\max}$ at the end of the collision
process we can imagine that the particles would be more
thermalized with respect to the ones that come from a tree where
$L(k)<m_{\max}$. This off course is not true for definition
(\ref{def:1}).

In table (\ref{tab:len}) we show the values of the length of the
collision trees reported in Figure (\ref{fig:b}) using the three
different definitions.

The implementation of such strategy inside the recursive algorithm
needs a modification of the first formulation given in section
\ref{sec:VHS} because it is necessary to evaluate the length of
collision trees without performing collisions. The idea is to
write into a list the collision process (performed in a recursive
way), to evaluate the length and to perform the collision using
the collision process stored into the list.

The following algorithm shows how we can store the collision tree
(correspondent to a process to sample from the density function
$f_k$) into a list named $tree$ and how we can evaluate his length
by using the recursive definition seen above.\\
\begin{algorithm}[Storage and length of collision trees]  ~\par
{\samepage\small
\begin{tabbing}
\ind 1. \= assign to a variable $indx$ the value $indx=0$  \\
\ind 2. \> put $path(indx)= k$ \\
\ind 3. \> {\fP if}  \= $k = 0$ return $0$ \\
                \> {\fP else} \= proceed as follows \\
                \>  \> - \= choose $j \in \left\{0,1, \ldots, k-1\right\}$ with equal probability \\
                \>  \> - \> put $h = k - j -1$ \\
                \>  \> - \> increment the index:  $indx = indx + 1$ \\
                \>  \> - \> store $j$ in the list: $path(indx) = j$ \\
                \>  \> - \> increment the index: $ indx = indx + 1$ \\
                \>  \> - \> store $h$ in the list: $path(indx) = h$ \\
                \>  \> - \> return $(1 + $ \\
                \>  \>   \> $+ \min\{($repeat from step 3 with $k = j), ($repeat from step 3 with $k =
                h)\}$\\
\end{tabbing}
\label{alg:2-r}
}
\end{algorithm}

In this way we obtain a list $path$ which contains the structure
of the collision tree and its relative length using definition
(\ref{def:2}). Off course the same algorithm applies also to
definition (\ref{def:3}).

Now we can finally modify the algorithm \ref{alg:1-r} at the point
5 in order to use the previously created collision trees.\\
\begin{algorithm}[TRMC-WB for VHS molecules] ~\par
{\samepage\small
\begin{tabbing}
\ind 1. \= compute the initial velocity of the particles, $ \left\{ v_i^0, i=1,\ldots,N  \right\}$\\
                \> by sampling them from the initial velocity $f_0$\\
\ind 2. \> split particles into collision sets as in algorithm (\ref{alg:split})\\
\ind 3. \> set counters $c_n = 0$ for $n=1,\ldots,m+1$\\
\ind 4. \> compute an upper bound $\Sigma$ of $\sigma_{ij}$\\
\ind 5. \> {\fP for} $n = m,\ldots,1$\\
            \> \ind \= take $N_n$ samples from the distribution with density $f_n$, according to\\
        \>          \> {\fP repeat}\\
        \>      \>  \ind - \= compute the list $path$ and its length $L$ as in algorithm \ref{alg:2-r}\\
        \>      \>  \ind - \> if $L>m_{\max}$ \\
        \>      \>         \> \ind \= set $N_{m+1}= N_{m+1}+1$ and $N_n=N_n - 1$\\
        \>      \>         \> else \\
        \>      \>         \>      \> set $indx = 0, k = path(indx) $\\
        \>      \>         \>   \ind 5.1. \= if $k=0$ sample $v_i$ from the initial density $f_0$ \\
        \>      \>         \>            \> else \= $indx=indx+1$ \\
        \>      \>         \>            \>      \> set $j = path(indx)$ \\
        \>      \>         \>            \>      \> set $indx = indx +1$ \\
        \>      \>         \>            \>      \> set $h = path(indx)$ \\
        \>      \>         \>            \>      \> if $c_j>0$ \= use a stored particle with a random choice\\
        \>      \>         \>            \>      \>            \> set $c_j = c_j - 1$, $N_j=N_j+1$ and $indx=indx+2j$\\
        \>      \>         \>            \>      \> else       \> sample $v_i$ and $v_i^*$ from $f_j$ (recursively from 5.1\\
        \>        \>       \>            \>      \>            \>  by setting $k=j$)\\
        \>        \>       \>            \>      \>            \>  $v_i^*$ is stored and then set $c_j=c_j+1$, $N_j=N_j-1$\\
        \>      \>         \>            \>      \> end if\\
        \>        \>       \>            \>      \> if $c_h>0$ \> use a stored particle with a random choice\\
        \>        \>       \>            \>      \>            \> set $c_h = c_h - 1$, $N_h=N_h+1$ and $indx=indx+2h$\\
        \>        \>       \>            \>      \> else       \> sample $v_j$ and $v_j^*$ from $f_h$ (recursively from 5.1\\
        \>        \>       \>            \>      \>             \>  by setting $k=h$)\\
        \>        \>       \>            \>      \>             \>  $v_j^*$ is stored and then set $c_j=c_j+1$, $N_j=N_j-1$\\
        \>      \>         \>            \>      \> end if\\
        \>        \>       \>            \>      \>  Perform the dummy collision between $v_i$ and $v_j$ as in
        DSMC.\\
        \>      \>         \>            \> end if\\
        \>        \>       \> \>set $N_n=N_n - 2$ \\
        \>        \>       \> end if\\
        \>        \>      {\fP until} $(N_n>0)$ \\
\ind    \> {\fP end for} \\
\ind 6. \> sample $N_{m+1}$ particles from the local Maxwellian\\

\end{tabbing}
\label{alg:3-r}
}
\end{algorithm}

\begin{figure}[t]
\centering
\includegraphics[scale=0.50]{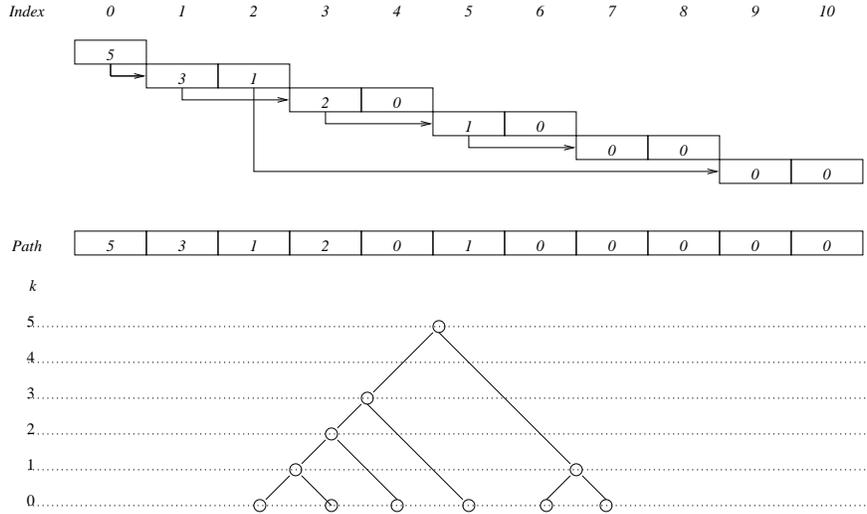}
\caption{Example of storage of a collision tree}
\label{fig:stor_tree}
\end{figure}

\begin{remark}
In the Variable Hard Sphere case the calculation of algorithm
\ref{alg:2-r} does not represent the effective length of the real
collision trees, due to the dummy collisions which can occur
during the total collision process. Thus the collisional trees
performed in VHS model are usually longer than the ones for
Maxwell molecules.

\end{remark}

\section{Numerical tests}

We present some numerical tests, both in space homogeneous and
space non homogeneous situations using the Hard Sphere model. The
solution obtained by the Recursive Time Relaxed scheme combined
with the adaptive strategy, is compared with the one obtained by
the classic Bird's algorithm.

To simplify notation we will use TRMC-R for the basic Recursive
Time Relaxed Monte Carlo Scheme defined by Algorithm 4.4, TRMC-RAD
for same scheme improved by the adaptivity strategy of section 4.3
(based on the shear stress tensor as equilibrium indicator) and
TRMC-WB for the Recursive Time Relaxed Monte Carlo Scheme defined
by Algorithm 4.6 (based on the well balanced truncation of the
collision trees). With BIRD we refer to the classic DSMC scheme
\cite{Bird3-r}.

All macroscopic quantities have been considered in non dimensional
form. As efficiency indicator we have considered the total number
of collisions performed in the simulation. Sampling two particles
from the Maxwellian has been counted as one collision. Note that
due to the structure of the recursive algorithm in principle we
may have additional computational requirements in terms of storage
of particles during the evaluation of collision trees and storage
and calculation of the collision trees if we want to implement a
shape-depending truncation of the trees. However particle storage
along the collision trees does not require any additional memory.
In fact to achieve conservations each particle is used only one
time in a collision tree and so it is the particle itself which is
stored in the list. On the other hand the additional use of memory
for the storage of the computational trees involves only one
single collision history and thus the impact in the overall
calculation is negligible.

\subsection{Space homogeneous tests}

We consider the sum of two Maxwellian as initial data. The
solutions has been obtained in one single run using $5 \times
10^4$ particles and choosing a collision time step $\Delta
t/\epsilon$ equal to one. The reference solution, called in the
plot as "exact" has been obtained by a large number of averages on
Bird's scheme.

\begin{figure}[t]
\centering
\includegraphics[scale=0.37]{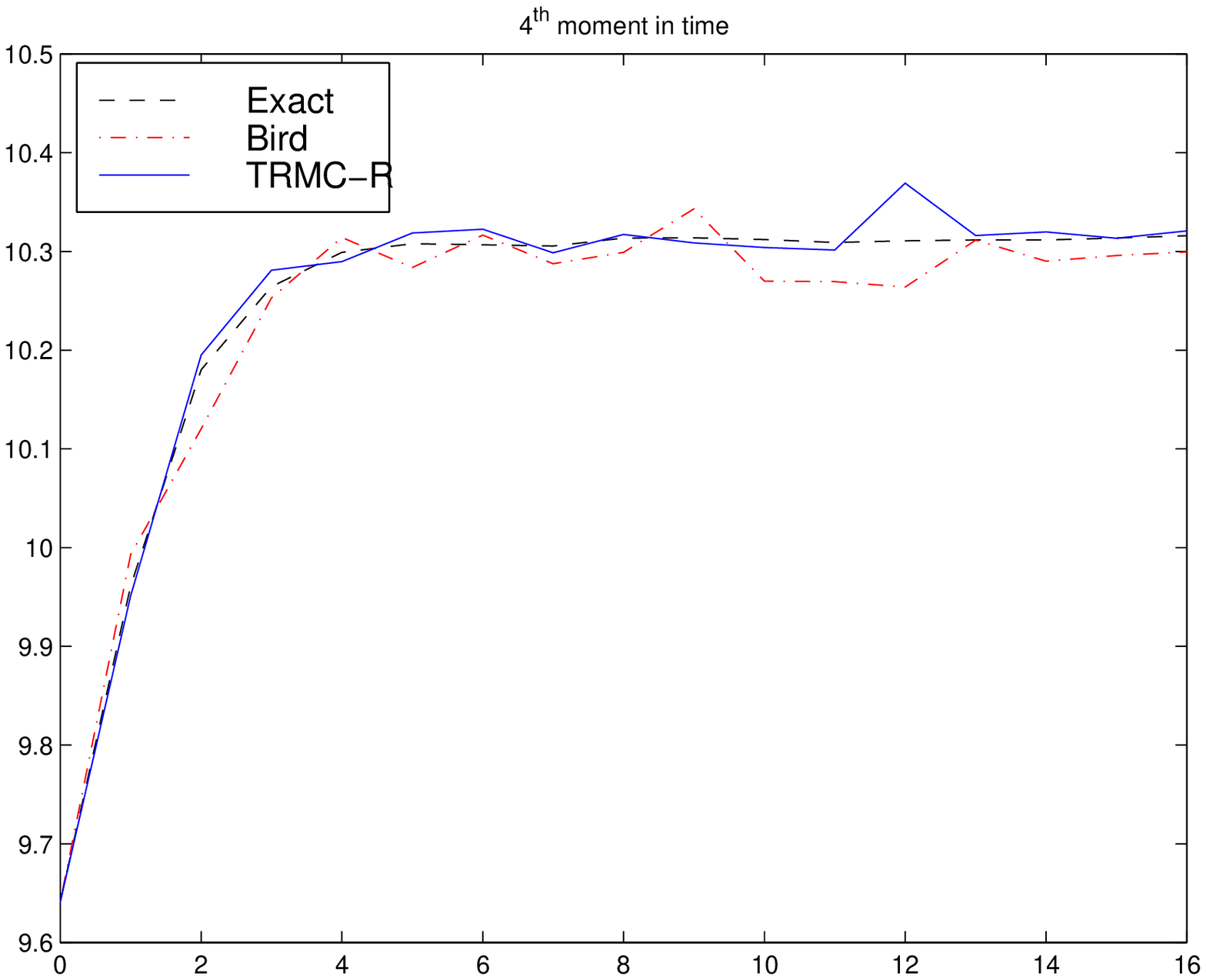}
\includegraphics[scale=0.37]{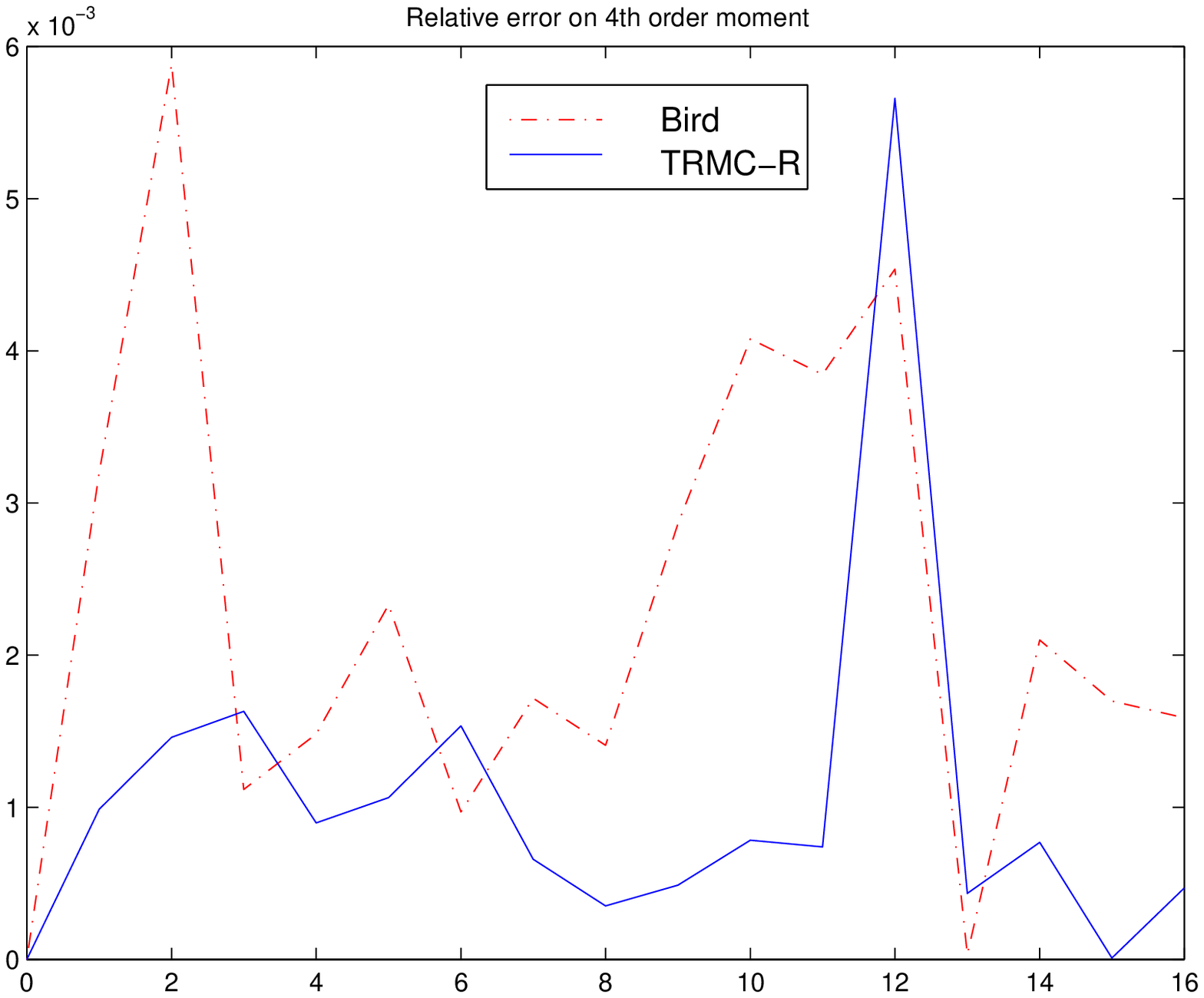}
\caption{Evolution in Time of $4^{th}$ order moment (left)  and
relative error (right) for TRMC-R and Bird's method.}
\label{fig:1}
\end{figure}

We compare the results for the fourth order moment $M_4$ and the
component $P_{xx}$ of the shear stress tensor
\[
M_4=\intrtre f v^4\,dv,\qquad P_{xx}=\intrtre f (v_1-u_1)^2\,dv.
\]

\begin{figure}
\centering
\includegraphics[scale=0.37]{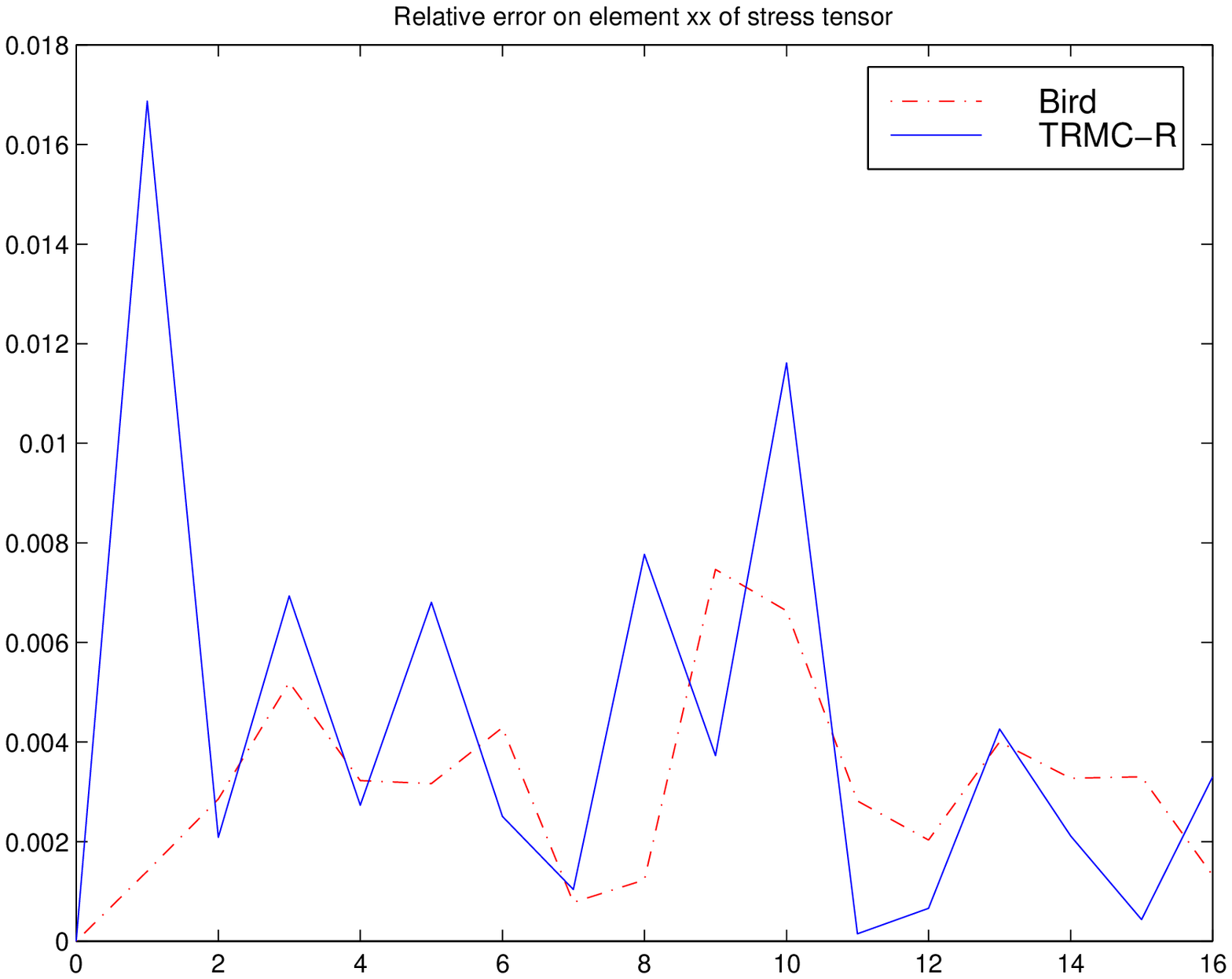}
\includegraphics[scale=0.37]{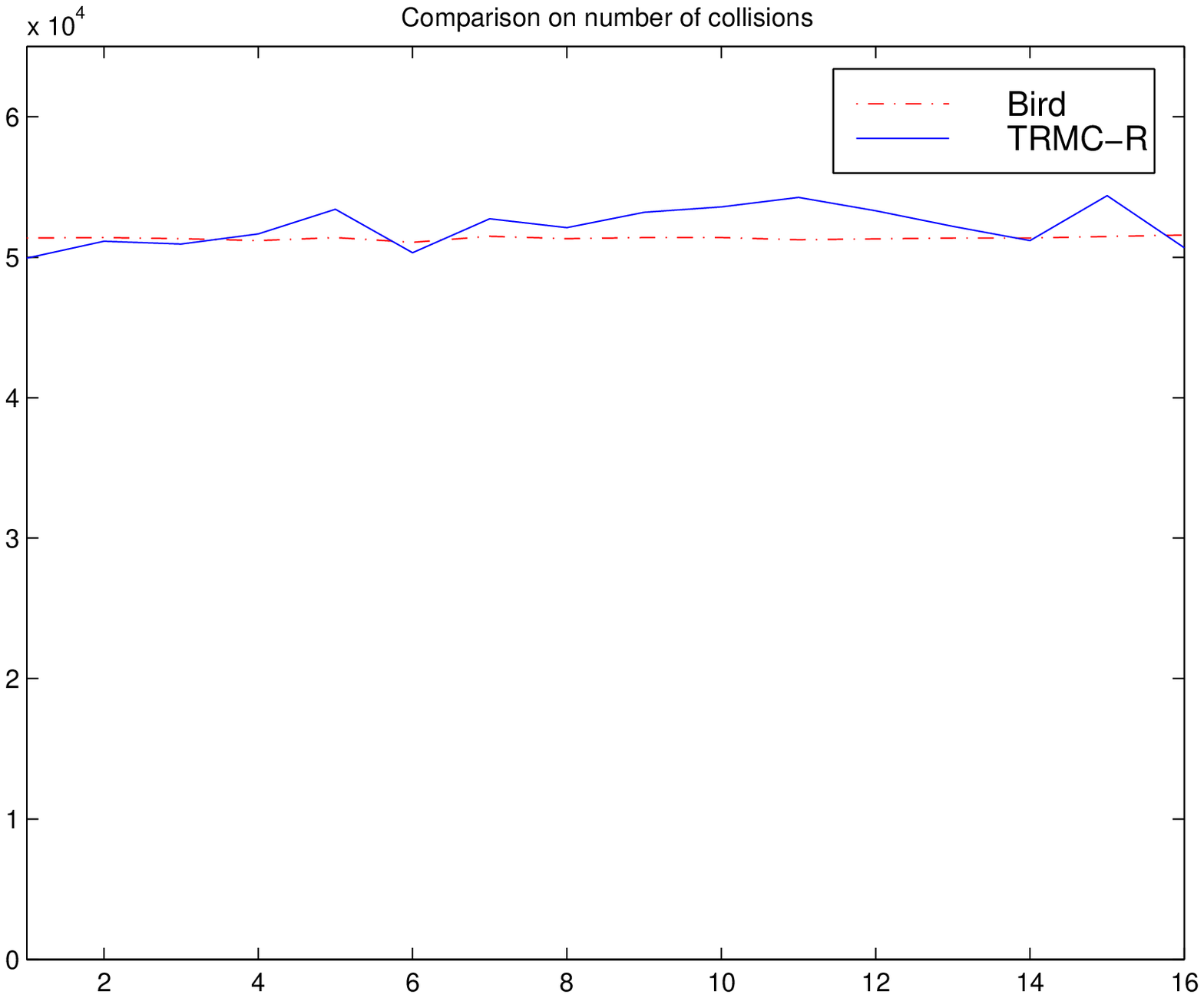}
\caption{Relative error for the evolution in time of the stress
tensor component $P_{xx}$ (left) and number of collisions (right)
for TRMC-R and Bird's method.} \label{fig:2}
\end{figure}

Figure \ref{fig:1} shows the case of TRMC-R without any limit for
the maximum depth of collision trees. The agreement with respect
to Bird's scheme is very good, and if we look at the relative
errors we can observe that they are comparable. In Figure
\ref{fig:2} we plot the relative error for the shear stress tensor
component $P_{xx}$ and the number of collisions necessary to
perform the simulations. We can consider the number of collisions
as an index of efficiency in terms of computational time. A single
sample from a Maxwellian is counted as half of a collision. As
expected we had no gain in efficiency because with this scheme we
took into account all possible collision processes.

\begin{figure}[t]
\centering
\includegraphics[scale=0.37]{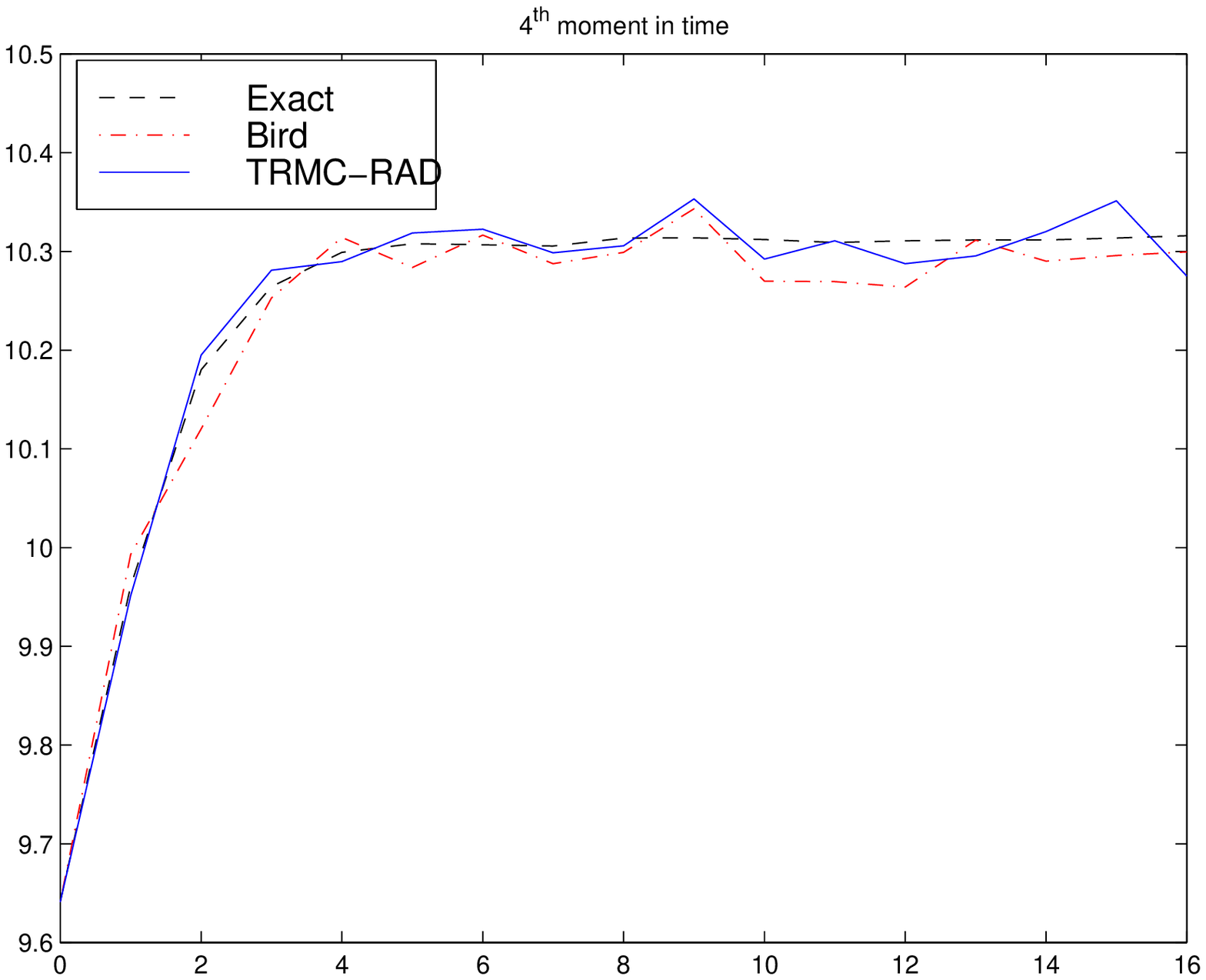}
\includegraphics[scale=0.37]{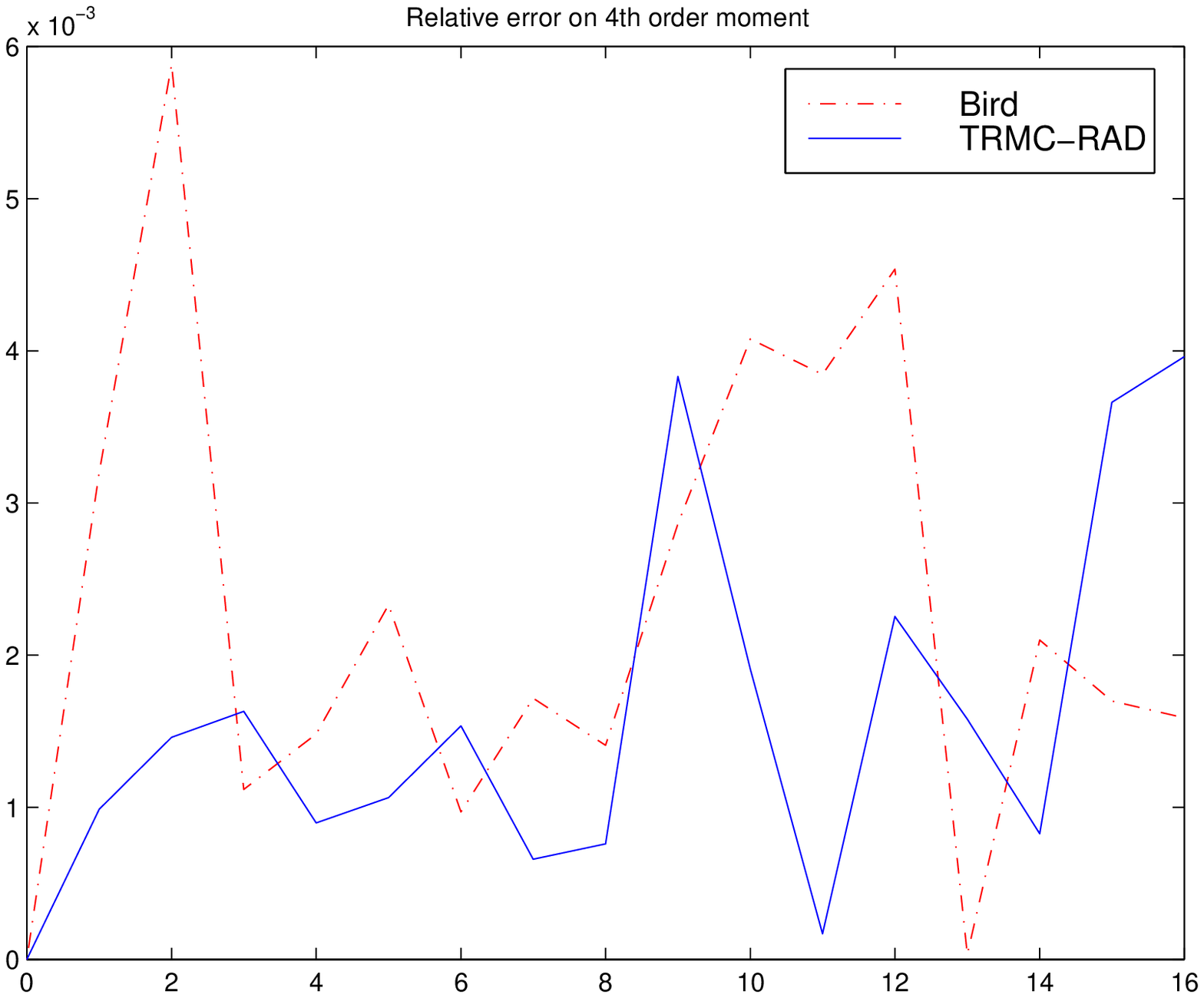}
\caption{Evolution in Time of $4^{th}$ order moment (left)  and
relative error (right) for TRMC-RAD and Bird's method.}
\label{fig:3}
\end{figure}

\begin{figure}
\centering
\includegraphics[scale=0.37]{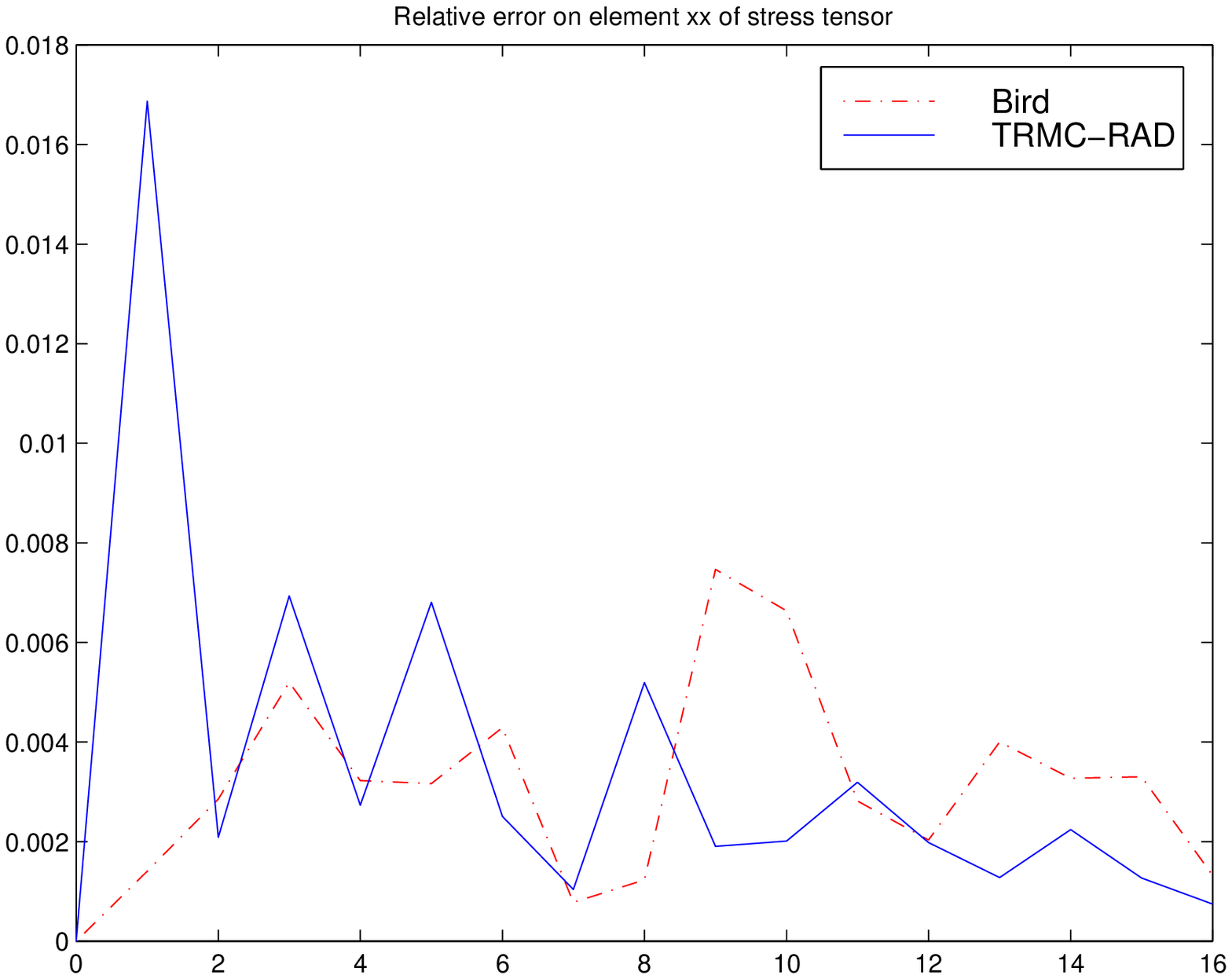}
\includegraphics[scale=0.37]{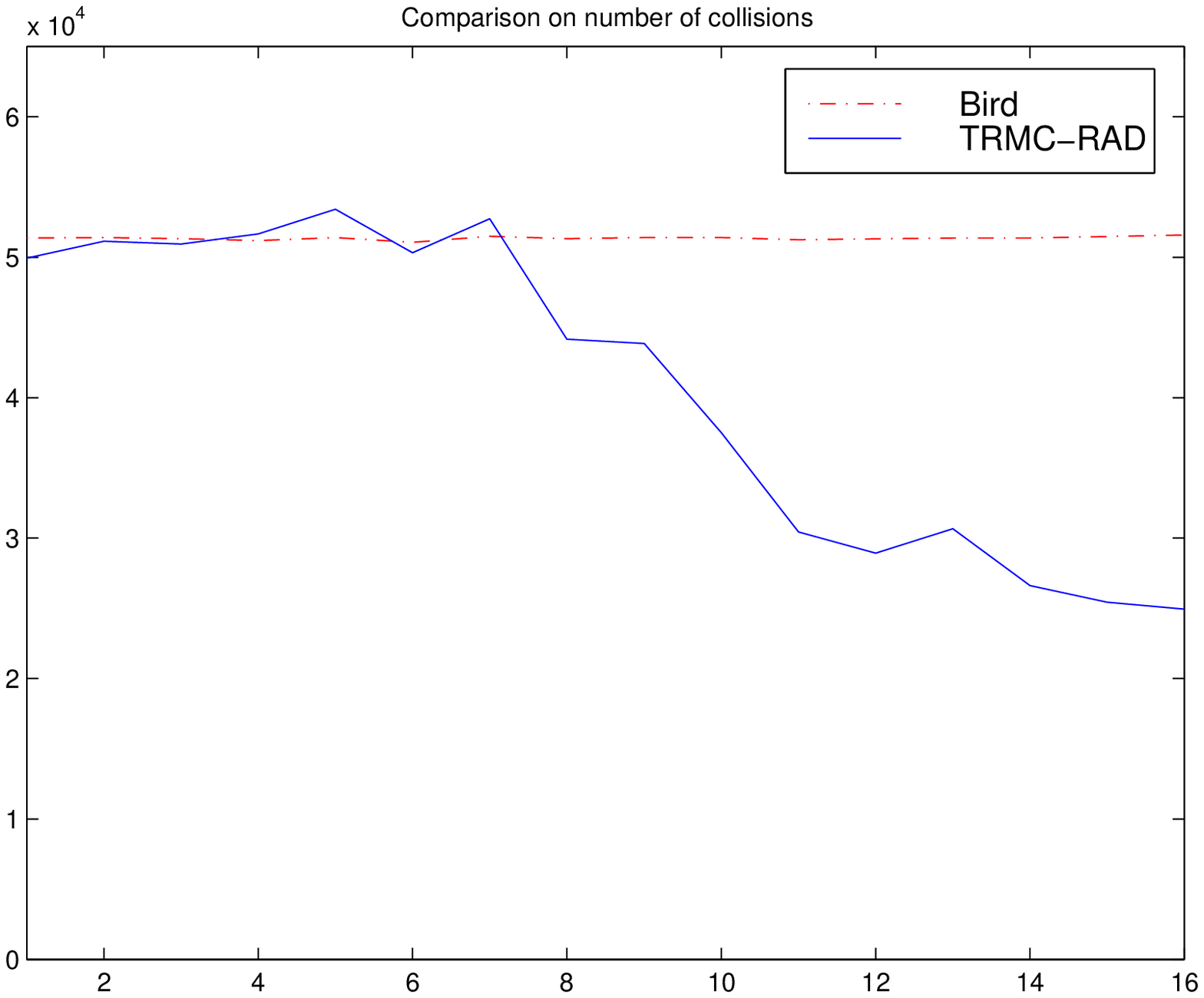}
\caption{Relative error for the evolution in time of the stress
tensor component $P_{xx}$ (left) and number of collisions (right)
for TRMC-RAD and Bird's method.} \label{fig:4}
\end{figure}

In the next case (Figure \ref{fig:3} and \ref{fig:4}) we perform
the same simulation using TRMC-RAD, starting with $m_{max} = 2$
and taking $\delta_1 = 0.005$ and $\delta_2 = 0.01$. We preserve
accuracy with respect to DSCM solution, but we have obtained a
strong reduction of the computational time because a bigger
fraction of particles is sampled directly from the Maxwellian
without going through the whole collision tree.

\begin{figure}[t]
\centering
\includegraphics[scale=0.37]{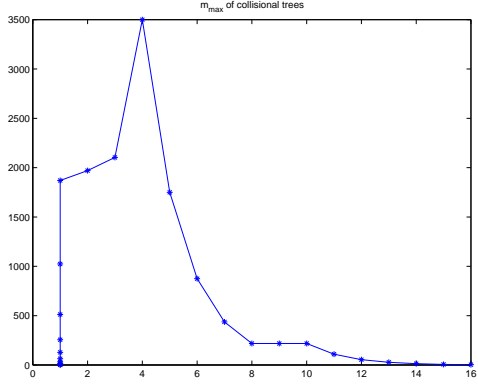}
\caption{Maximum depth $m_{max}$ of a collision tree in time for
TRMC-RAD} \label{fig:5}
\end{figure}

We observe in Figure \ref{fig:5}, where we report the maximum
depth of the collision trees during the calculation, that the
simulation in the first time steps has been performed several
times, using at each time $m_{\max} = 2 m_{\max}$. As already
mentioned let us point out that these calculations do not affect
the total computational cost because all ``old collisions" are
kept and reused in order to complete the new collision process.

\begin{figure}[t]
\centering
\includegraphics[scale=0.37]{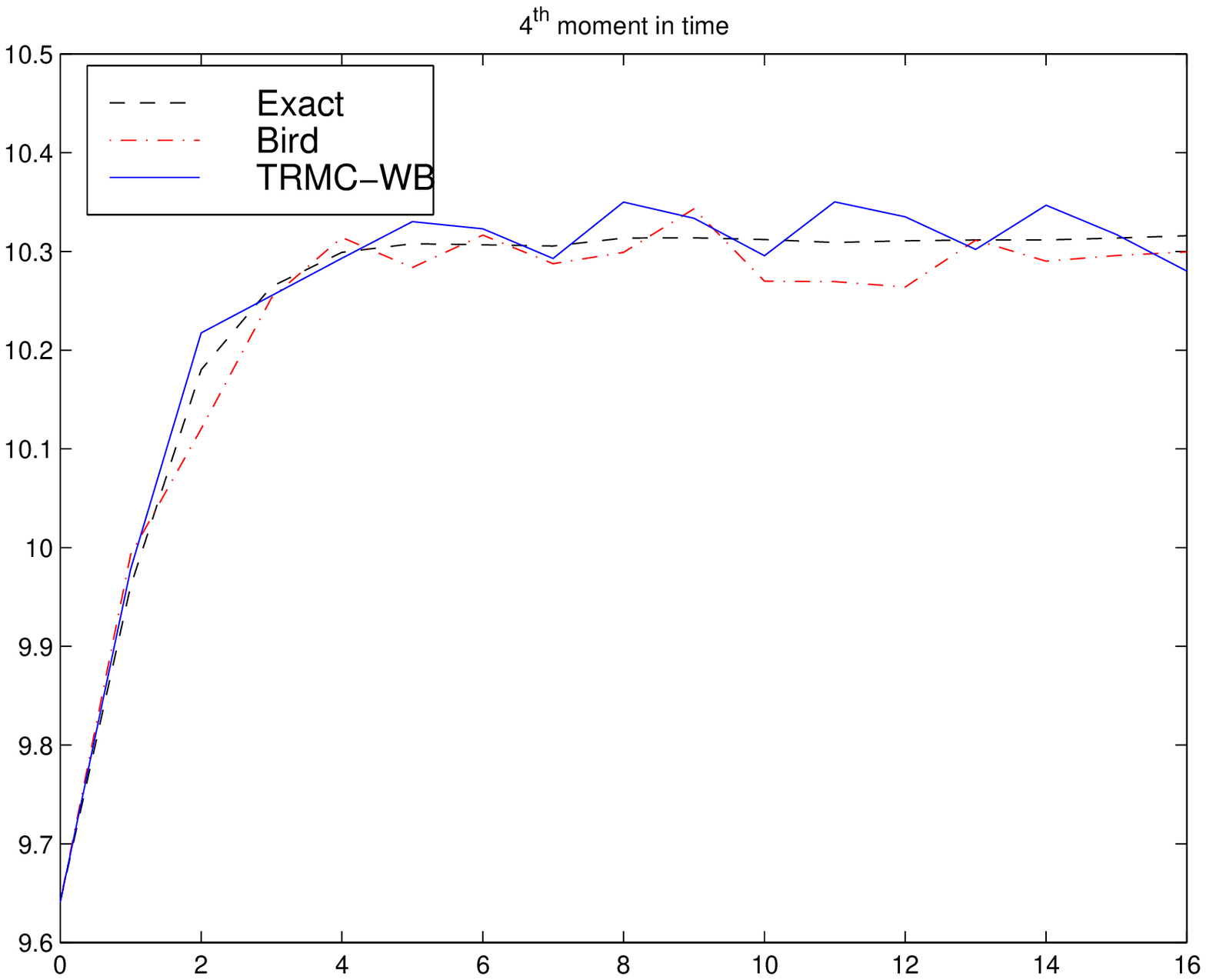}
\includegraphics[scale=0.37]{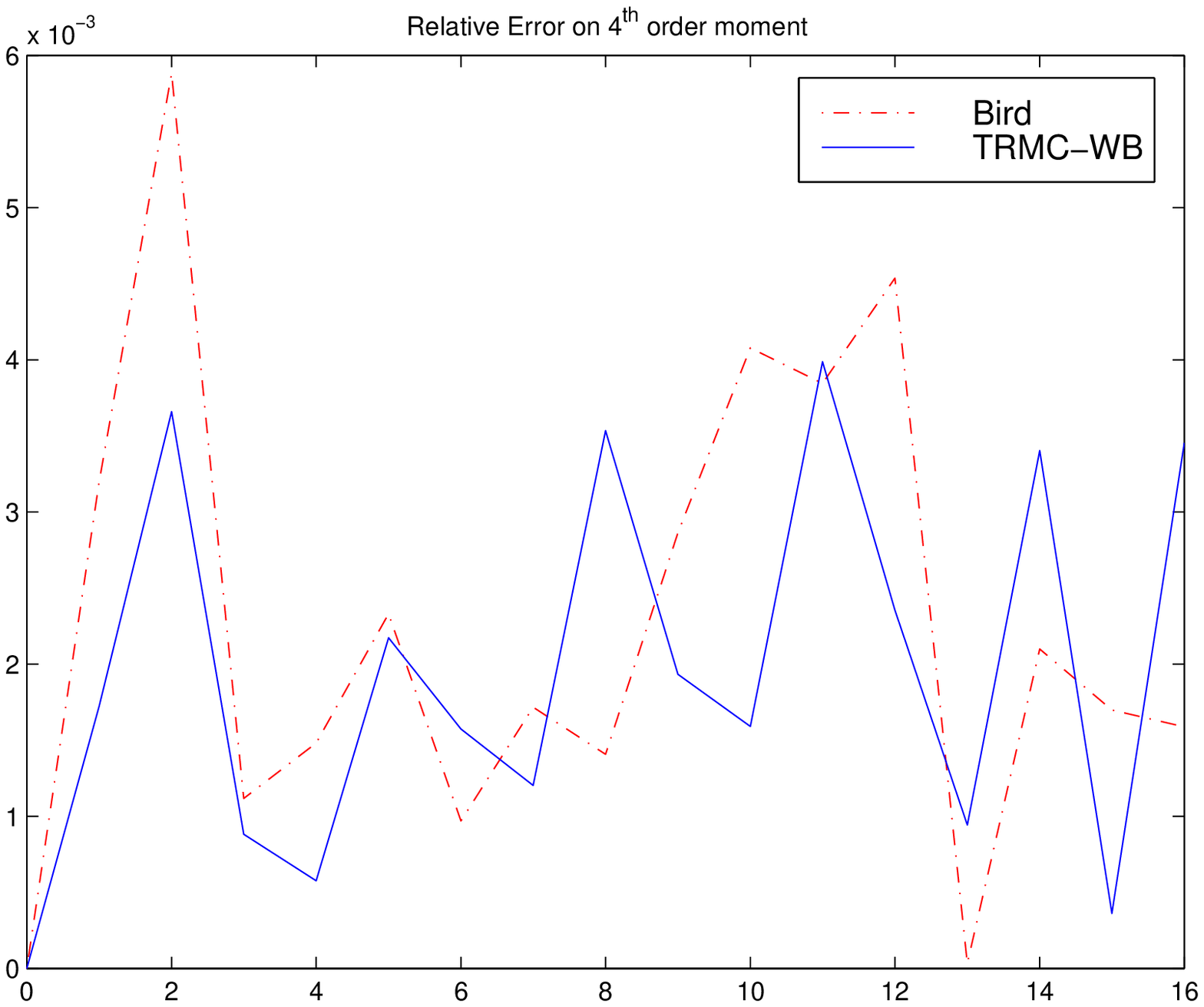}
\caption{Evolution in Time of $4^{th}$ order moment (left)  and
relative error (right) for TRMC-WB and Bird's method.}
\label{fig:6}
\end{figure}

\begin{figure}
\centering
\includegraphics[scale=0.36]{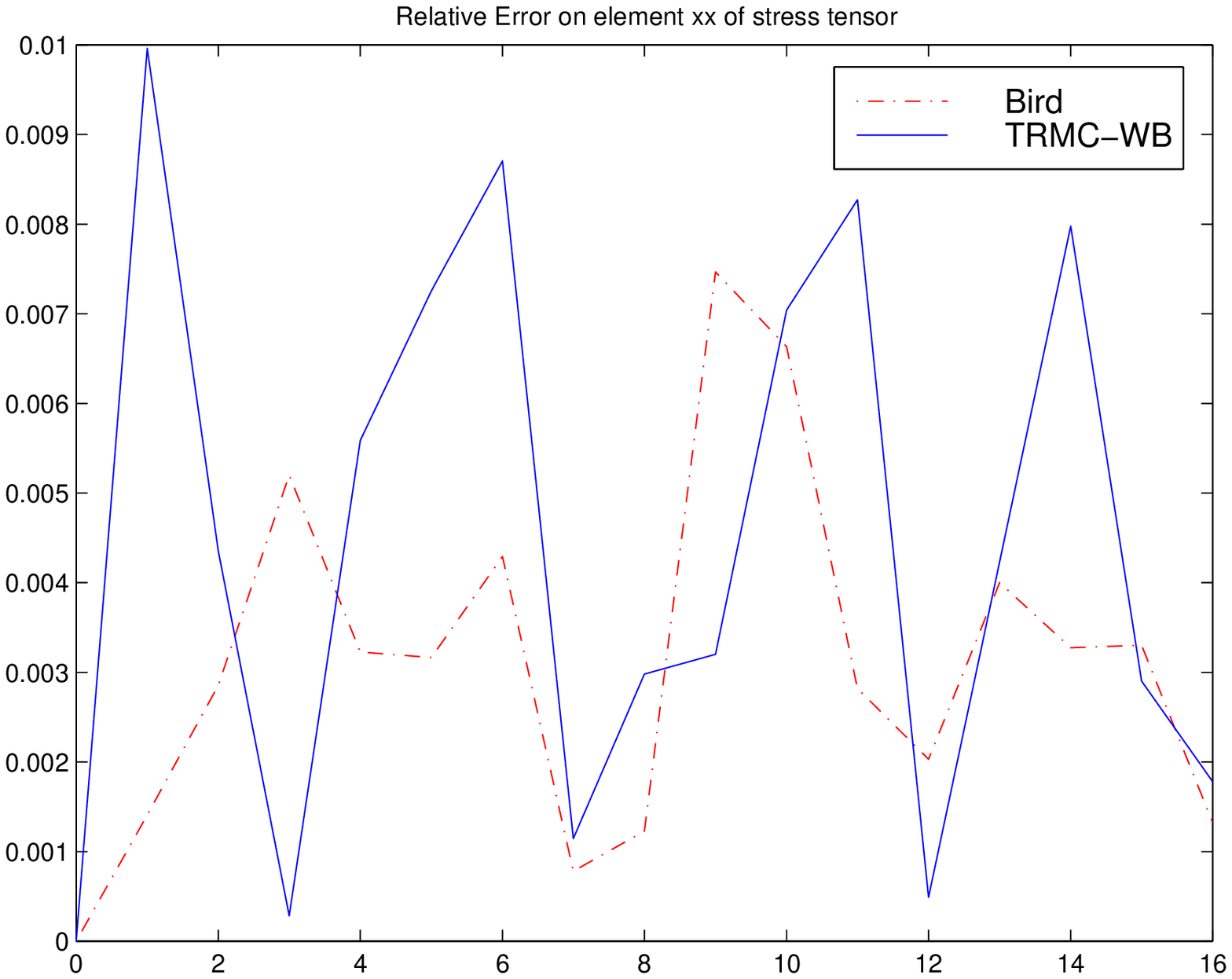}
\includegraphics[scale=0.36]{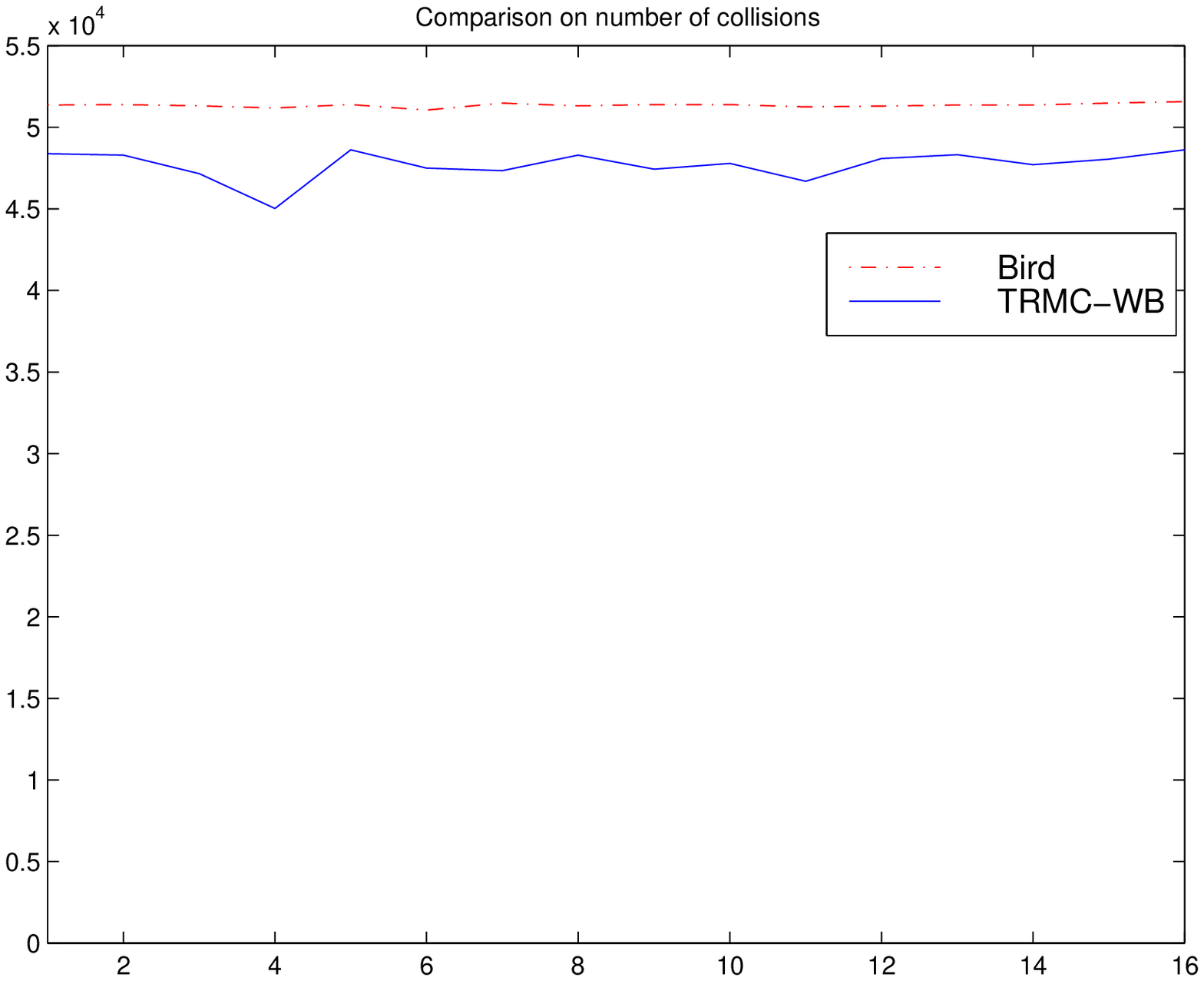}
\caption{Relative error for the evolution in time of stress tensor
component $P_{xx}$ (left) and number of collisions (right) for
TRMC-WB and Bird's method.} \label{fig:7}
\end{figure}

The last homogeneous test deals with the shape depending strategy
based on truncation of collision trees by using the definition
$L(k = h+j+1) = 1 + \min\{L(h),L(j)\}$  and $m_{\max} = 5$.
Similar results can be obtained using the definition $L(k = h+j+1)
= 1 + {\rm mean}\{L(h),L(j)\}$, we omit them for brevity. We have,
also in this case, a gain in efficiency, preserving the accuracy
in time (see Figures \ref{fig:6} and \ref{fig:7}).

The TRMC-WB scheme shows a constant gain of computational cost
with respect to Bird's scheme, approximatively the 7\% during the
whole simulation, while the TRMC-RAD achieves the maximum gain
(about 50\%) at the end of the calculation, with an average gain
of 20\% in the total simulation time. The combination of TRMC-WB
and TRMC-RAD is under investigation and will be presented
elsewhere.

\subsection{Stationary shock}

Next we consider a non homogeneous stationary shock problem. The
boundary conditions have been assigned at the left and right side
accordingly to the Rankine-Hugoniot relations

\[ \rho_L u_L = \rho_R u_R, \]
\[ \rho_L u_L^2 + p_L = \rho_R u_R^2 + p_R, \]
\[ u_L(E_L + p_L) = u_R(E_R + p_R). \]

The values used in the simulation are
\begin{itemize}

    \item $M_{a_L}=3$ (Mach Number of incoming flux)
    \item $T_L = 1$ (Temperature of incoming flux)
    \item $u_{x_L}= -M\sqrt{\gamma T_L},\,u_{y_L}=0,\,u_{z_L}=0 $ $\gamma=5/3$ (Mean velocity of incoming flux)
    \item $\rho_L=1$ (Total mass)
  \item $\epsilon=1$, $\epsilon=0.1$ and $\epsilon=0.001$.

\end{itemize}

The numerical solution has been obtained using TRMC-RAD and Bird's
method with $50$ spatial cells and $1000$ particles in each cell.
Reference 'Exact' solution has been performed by Bird's standard
DSMC method using $50$ spatial cells and $3000$ particles in each
cell. A detailed analysis on the effect of the number of particles
per cell in TRMC methods has been performed in \cite{homo-r,
couette-r}.

\begin{figure}
\centering
\includegraphics[scale=0.37]{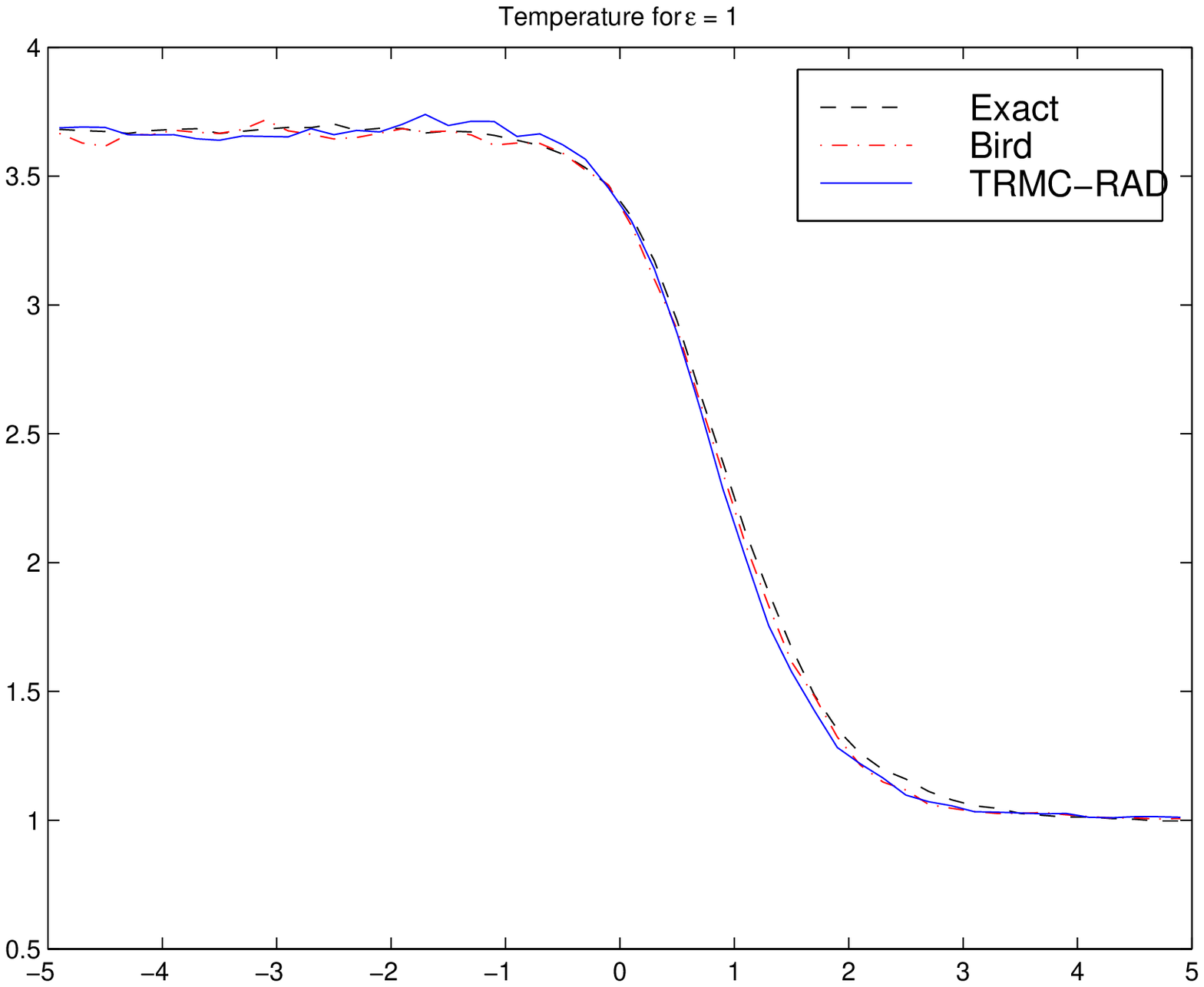}
\includegraphics[scale=0.37]{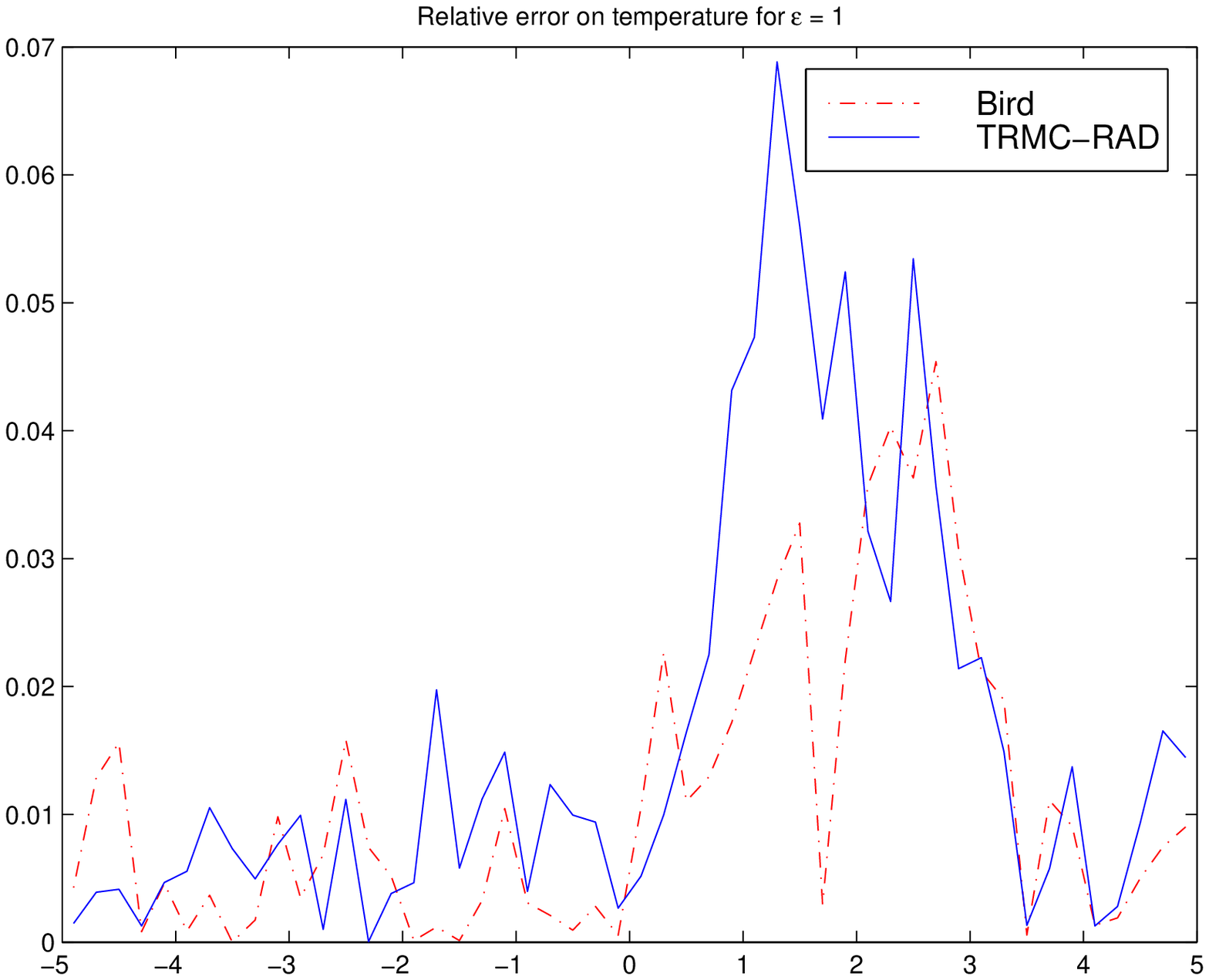}
\caption{Rarefied regime ($\epsilon=1$). Temperature (left) and
Relative Error (right) for TRMC-RAD and Bird's method.}
\label{fig:8}
\end{figure}

\begin{figure}
\centering
\includegraphics[scale=0.37]{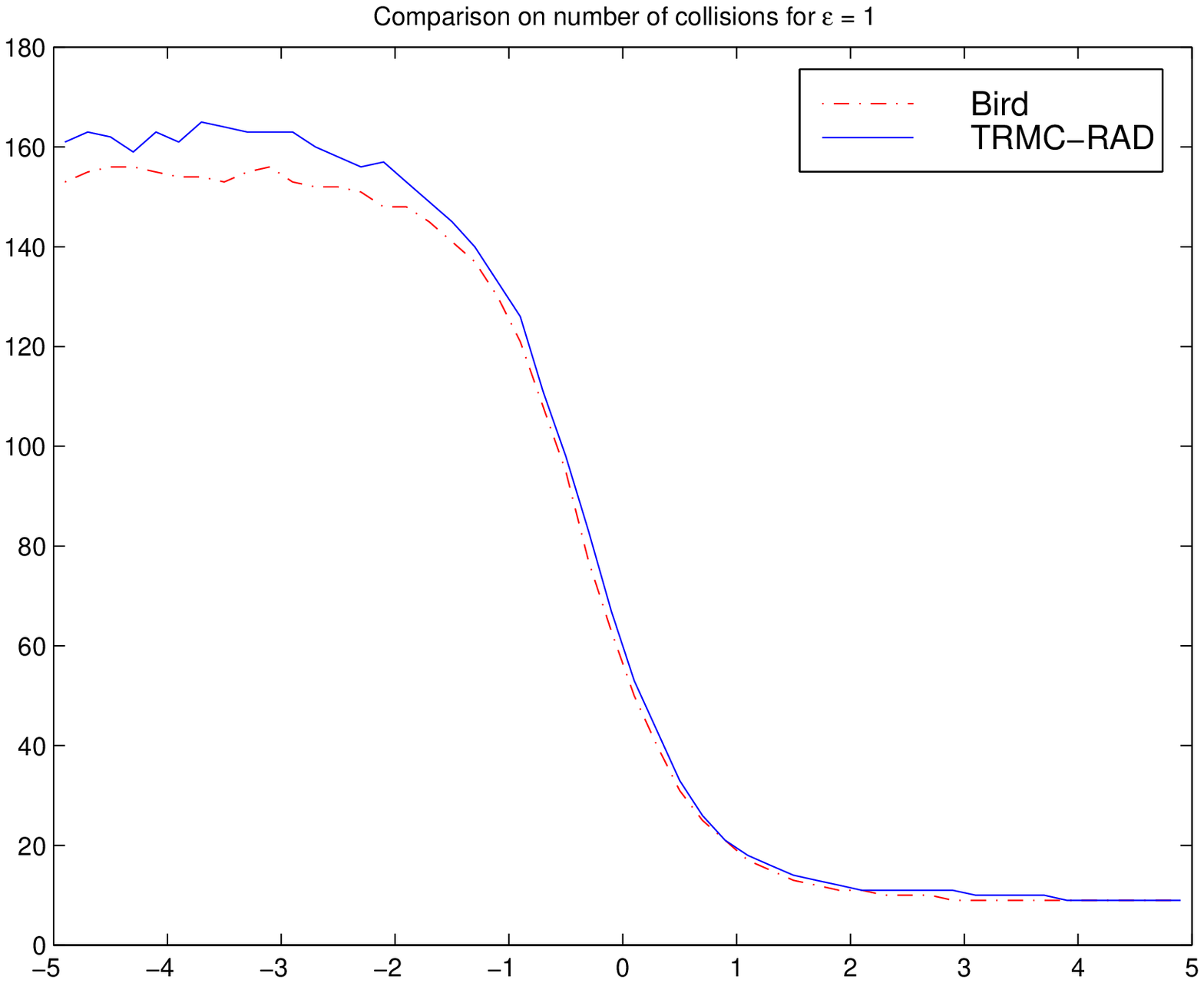}
\includegraphics[scale=0.37]{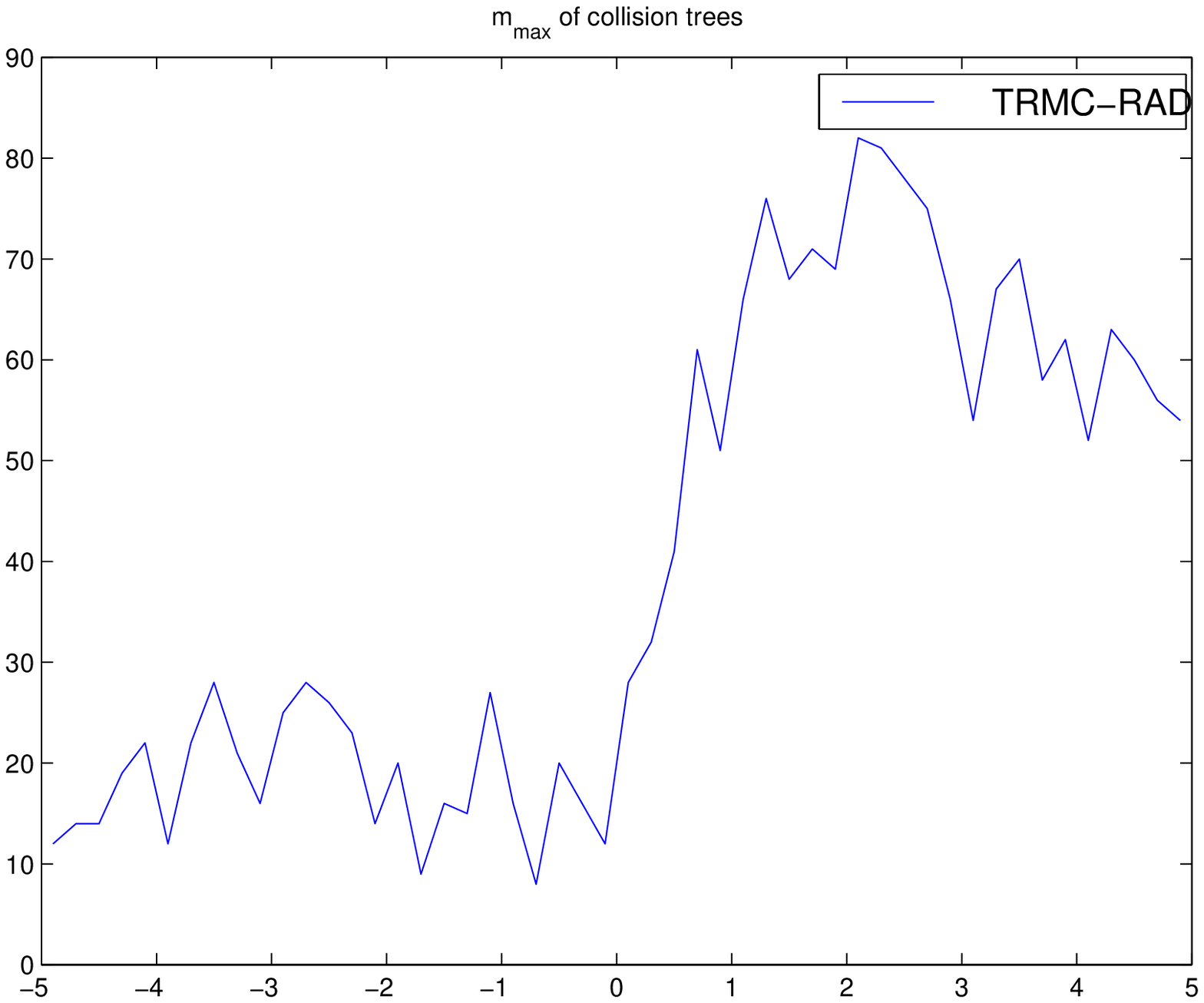}
\caption{Rarefied regime ($\epsilon=1$). Number of collisions for
TRMC-RAD and Bird's method (left) and maximum depth $m_{\max}$ in
each cell for TRMC-RAD (right). } \label{fig:9}
\end{figure}

\begin{figure}
\centering
\includegraphics[scale=0.37]{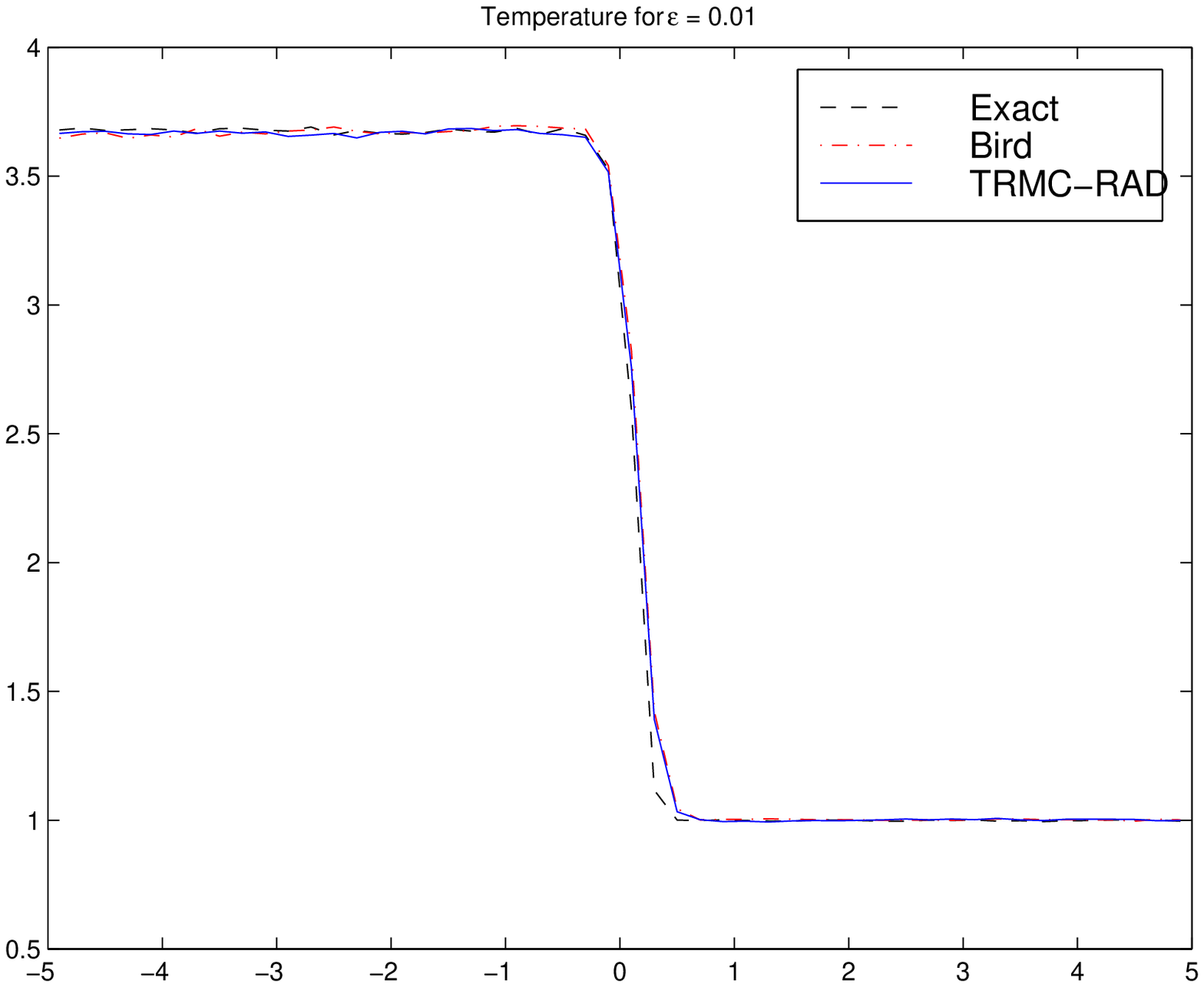}
\includegraphics[scale=0.37]{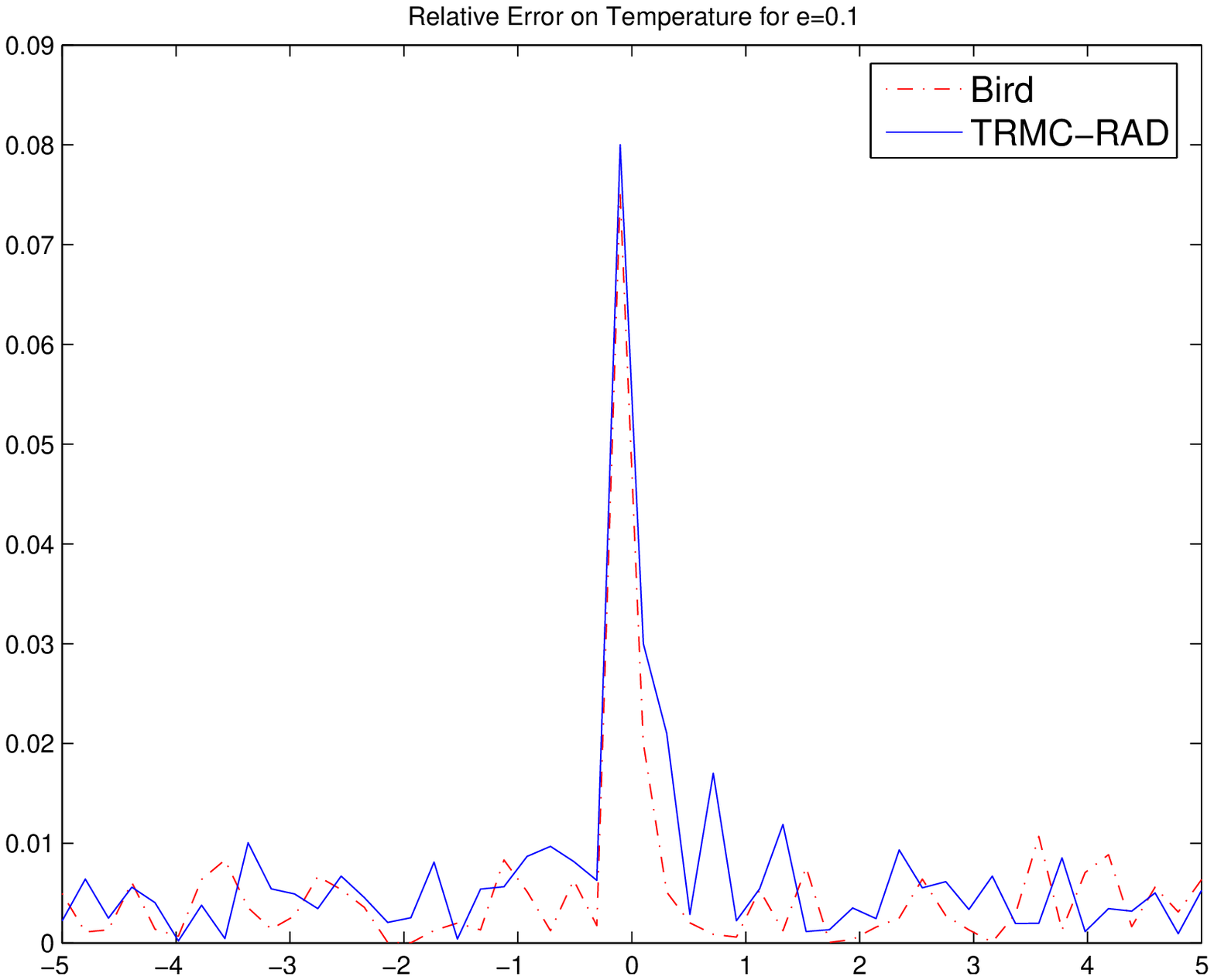}
\caption{Intermediate regime ($\epsilon=0.1$). Temperature (left)
and Relative Error (right) for TRMC-RAD and Bird's method.}
\label{fig:12}
\end{figure}

\begin{figure}
\centering
\includegraphics[scale=0.37]{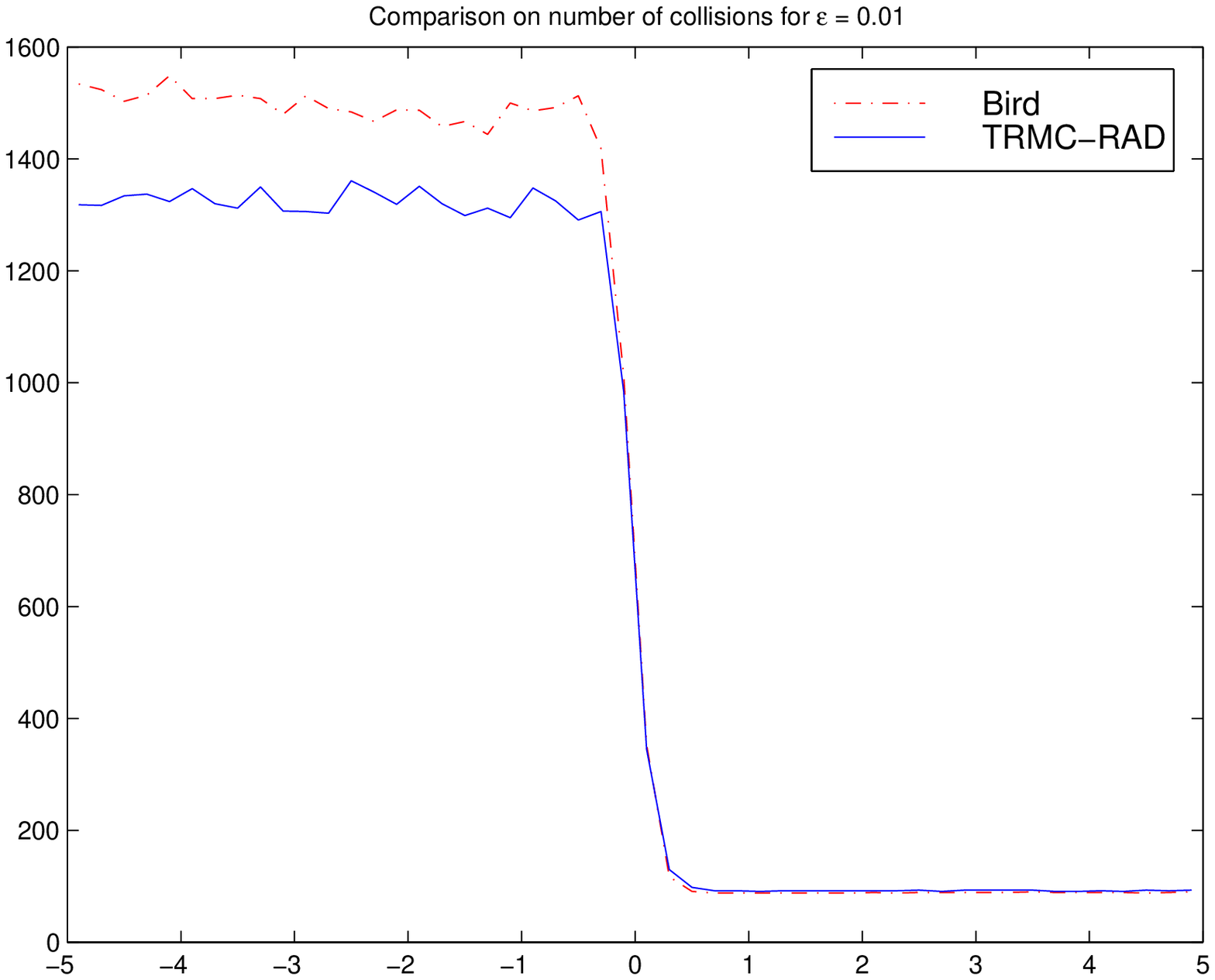}
\includegraphics[scale=0.37]{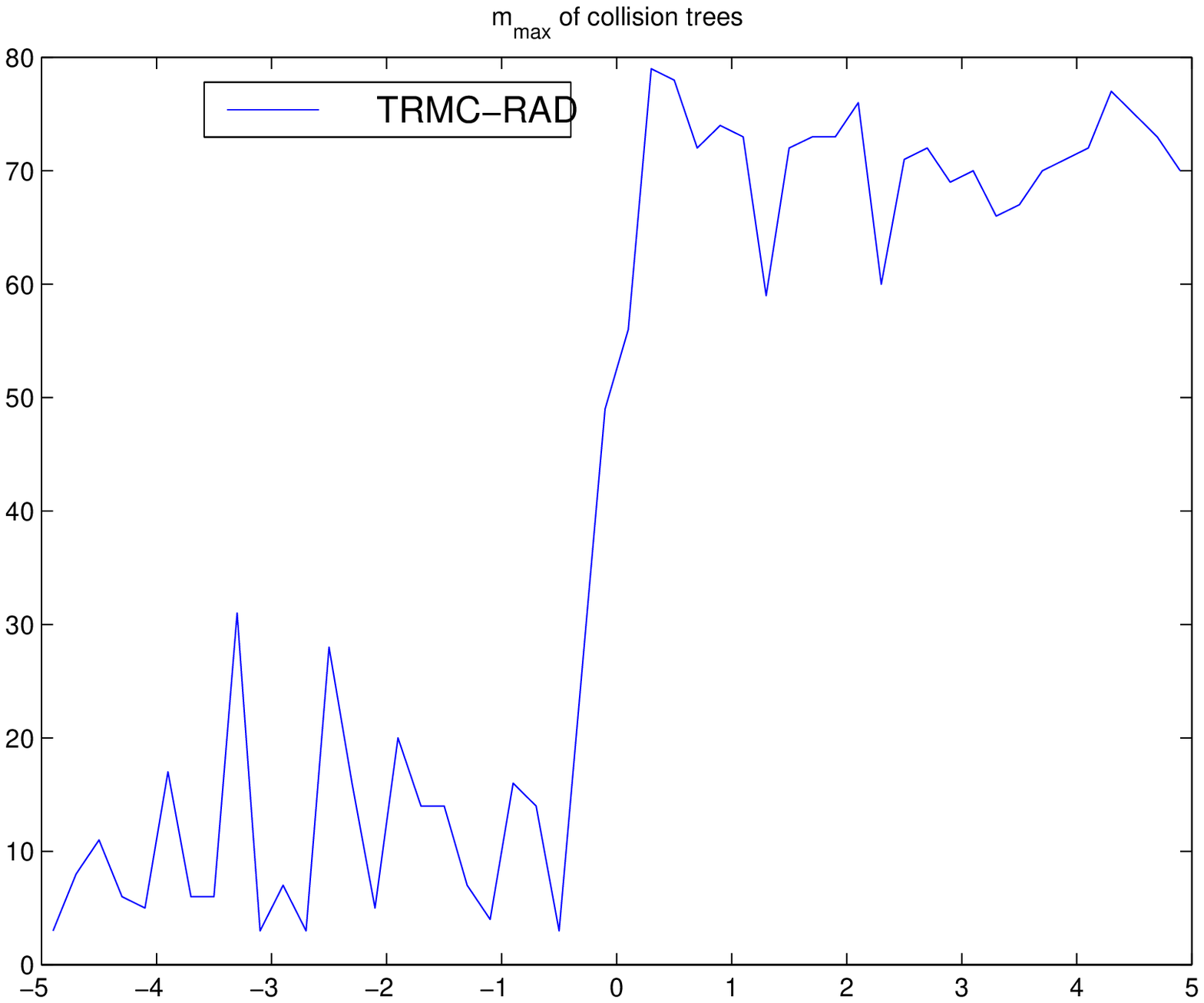}
\caption{Intermediate regime ($\epsilon=0.1$). Number of
collisions for TRMC-RAD and Bird's method (left) and maximum depth
$m_{\max}$ in each cell for TRMC-RAD (right). } \label{fig:13}
\end{figure}

\begin{figure}
\centering
\includegraphics[scale=0.37]{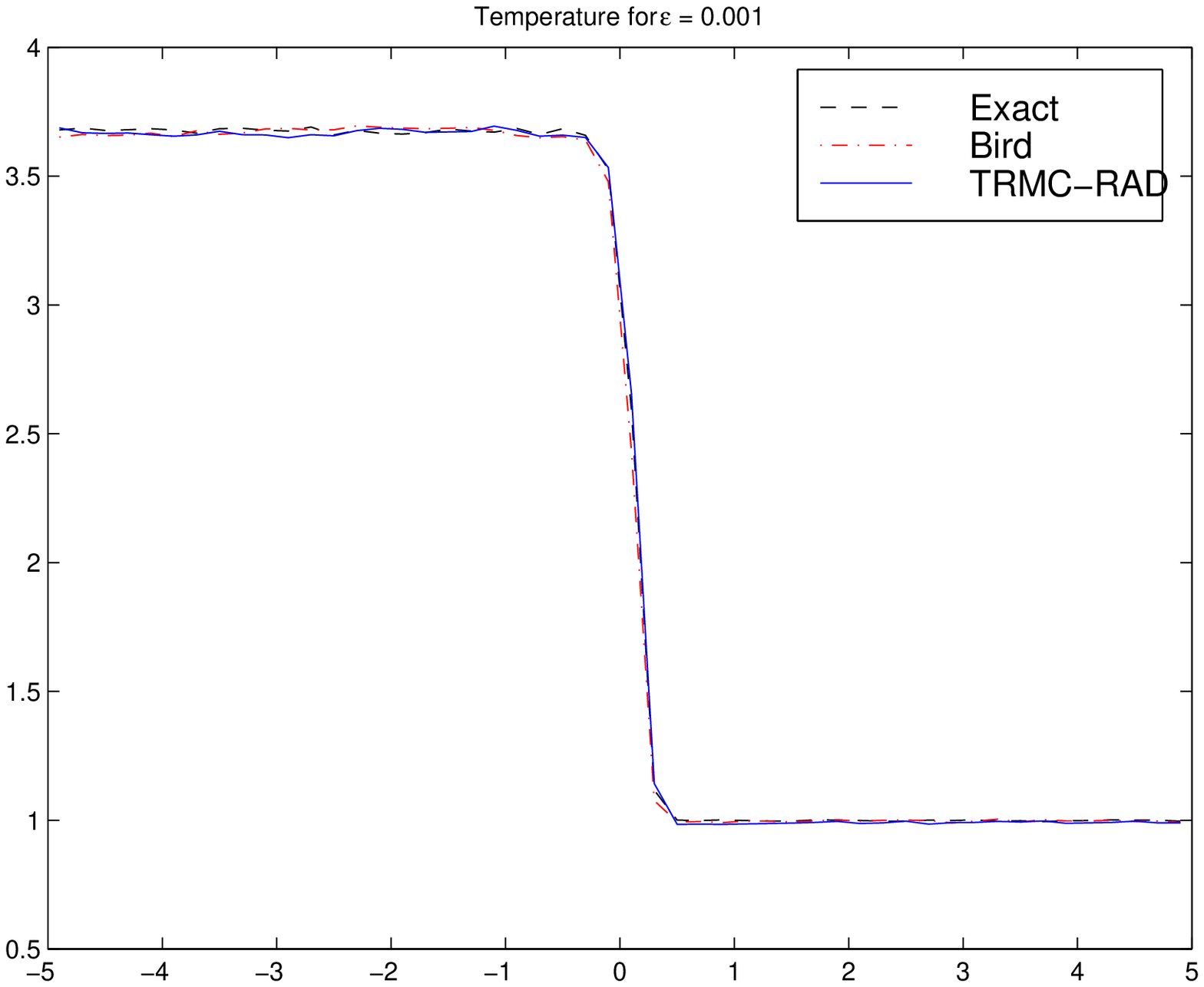}
\includegraphics[scale=0.37]{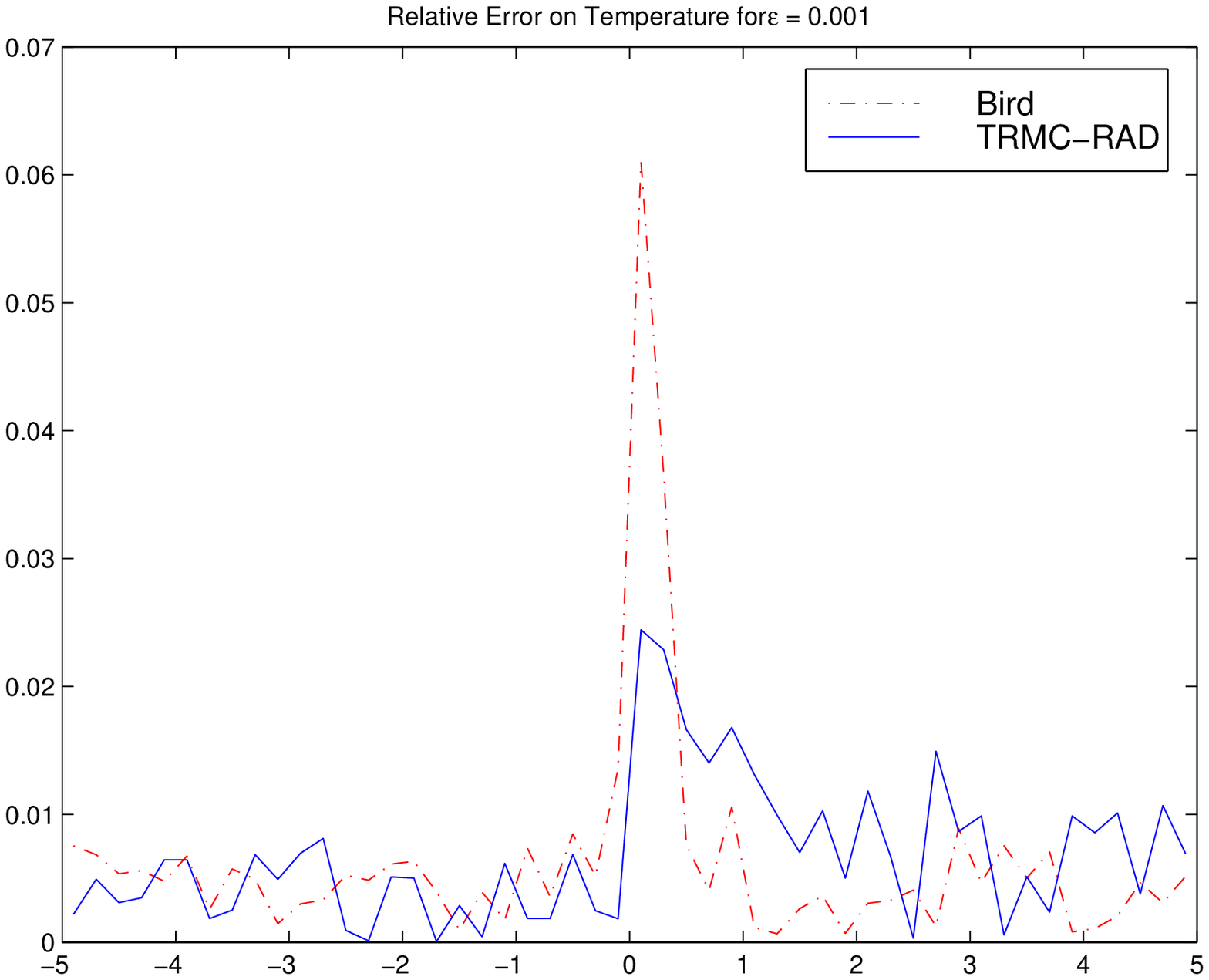}
\caption{Fluid regime ($\epsilon=0.001$). Temperature (left) and
Relative Error (right) for TRMC-RAD and Bird's method.}
\label{fig:10}
\end{figure}

\begin{figure}
\centering
\includegraphics[scale=0.37]{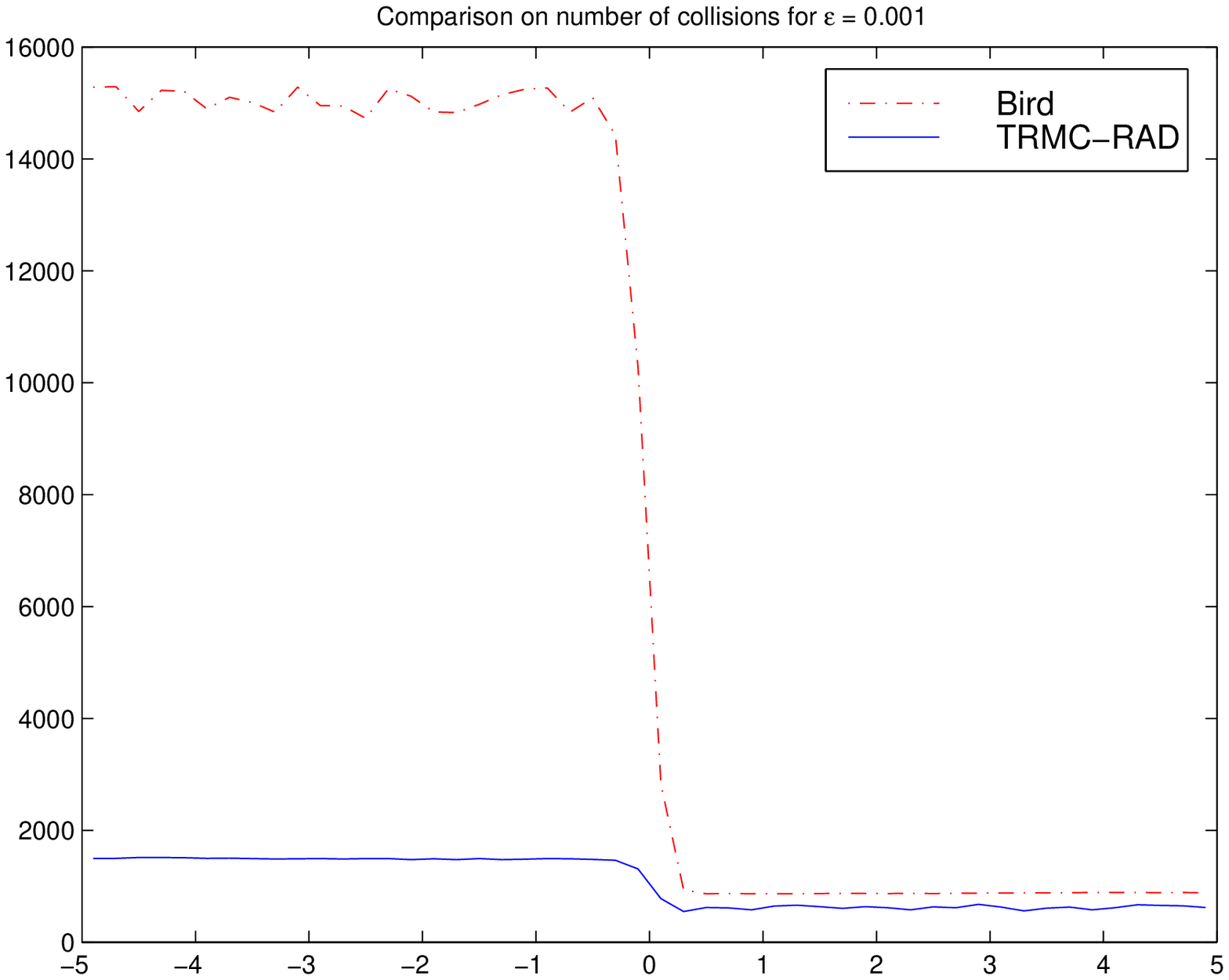}
\includegraphics[scale=0.37]{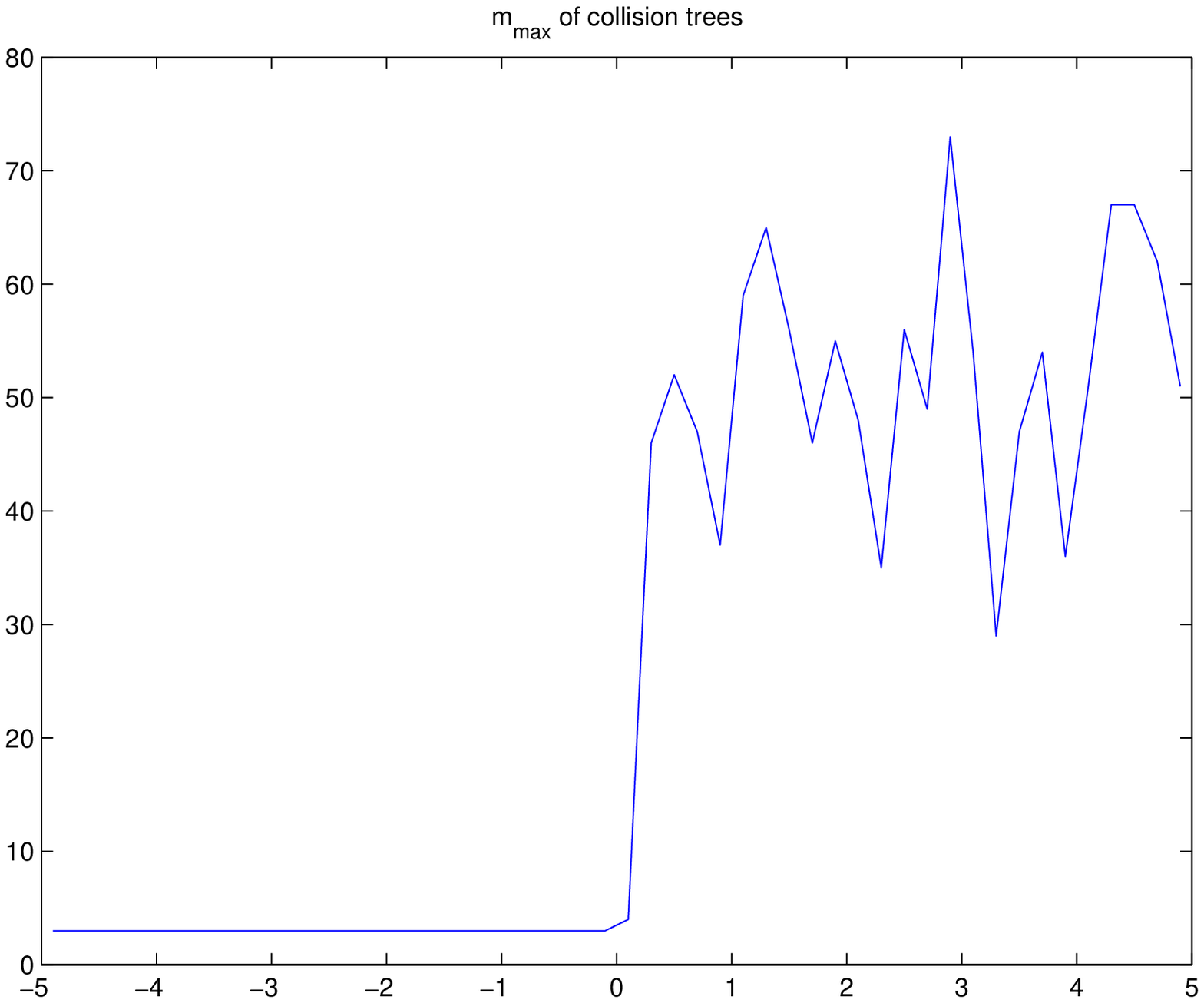}
\caption{Fluid regime ($\epsilon=0.001$). Number of collisions for
TRMC-RAD and Bird's method (left) and maximum depth $m_{\max}$ in
each cell for TRMC-RAD (right).} \label{fig:11}
\end{figure}

In order to increase accuracy, for $t$ large enough, averages of
the solution have been computed. The results for the temperature
shock are presented in Figures \ref{fig:8}-\ref{fig:11}. As
expected there is a good agreement between TRMC-RAD and Bird's
method and the relative errors are essentially comparable.

Note that there is an evident gain in efficiency of TRMC-RAD
against Bird's method without losing accuracy, especially near
fluid regime. Looking at the rarefied case $\epsilon =1$ we obtain
the same computational cost, while we have a gain close to 10\%
for the intermediate test case $\epsilon =0.1$ that increases up
to the 86\% close to the fluid regime for $\epsilon =0.001$.

\section{Conclusion}
We have presented recursive Monte Carlo methods which are suitable
for the numerical simulation of the Boltzmann equation for a wide
range of Knudsen numbers. These recursive TRMC methods minimize
the effects of  time discretization error and over-relaxation due
to the choice of the upper bound of the cross-section that were
present in the previous versions of TRMC. Of paramount importance
to increase the efficiency of the methods is the use of suitable
truncation strategies of the collisional trees, such as adaptive
truncation based on macroscopic quantities or well balanced
truncation based on the trees properties. The resulting schemes
are very promising and show that a considerable gain in efficiency
can be obtained without degradation in accuracy. However
additional test cases must be performed in order to validate the
schemes. A combination of the present schemes with the hybrid
strategy proposed in \cite{CPmc-r} is actually under
investigation.

% ----------------------------------------------------------------
\newpage
\bibliographystyle{amsplain}
%\bibliography{}

% ----------------------------------------------------------------
\end{document}